\documentclass{article}
\usepackage{amssymb}
\usepackage{amsmath}
\usepackage{amsfonts}
\usepackage{textcomp}
\usepackage{amsthm}
\usepackage{bookmark}
\usepackage{appendix}
\usepackage{dsfont}
\usepackage{changes}%\usepackage[final]{changes}  %\listofchanges
\usepackage{xcolor}
\usepackage{verbatim}
\usepackage{enumitem}
\usepackage{url}
\usepackage{hyperref}
\usepackage{cases}
\date{}

\newtheorem{Theorem}{Theorem}[section]
\newtheorem{Lemma}{Lemma}[section]
\newtheorem{Remark}{Remark}[section]
\newtheorem{Corollary}{Corollary}[section]

\newtheorem{Proposition}{Proposition}[section]

\numberwithin{equation}{section} \theoremstyle{plain}

\renewcommand{\theequation}{\arabic{section}.\arabic{equation}}

\allowdisplaybreaks[4]

\def\R{{\textbf{R}}}
\def\H{\mathcal{H}}

\def\f{\frac}
\def\R{\mathbb R}
\def\N{\mathbb N}
\def\C{\mathbb C}

\def\e{\mathrm{e}}

\def\al{\alpha}
\def\be{\beta}

\def\l{\Big}
\def\r{\Big}

\def\n{\nabla}
\def\g{\gamma}

\title{Optimizers of the Finite-Rank Hardy-Lieb-Thirring Inequality for Hardy-Schr\"odinger Operator}

\author{Bin Chen\thanks{School of Mathematics and Statistics,  Key Laboratory of Nonlinear Analysis $\&$ Applications (Ministry of Education), Central China Normal University, Wuhan 430079, P. R. China. Email:  binchenmath@mails.ccnu.edu.cn.}, Yujin Guo\thanks{School of Mathematics and Statistics,  Key Laboratory of Nonlinear Analysis $\&$ Applications (Ministry of Education), Central China Normal University, Wuhan 430079, P. R. China. Y. J. Guo is partially supported by National Key R $\&$ D Program of China (Grant 2023YFA1010001), and NSF of China (Grants 12225106 and 12371113). Email: yguo@ccnu.edu.cn.}, and Shuang Wu\thanks{School of Mathematics and Statistics,  and Hubei Key Laboratory of Mathematical Sciences, Central China Normal University, Wuhan 430079, P. R. China. Email: swu@mails.ccnu.edu.cn.}}

\date{\today}

\begin{document}
	\maketitle
	
\begin{abstract}
We study the following finite-rank Hardy-Lieb-Thirring inequality of Hardy-Schr\"odinger operator:
\begin{equation*}
			\sum_{i=1}^N\left|\lambda_i\Big(-\Delta-\f{c}{|x|^2}-V\Big)\right|^s\leq C_{s,d}^{(N)}\int_{\R^d}V_+^{s+\f d2}dx,
\end{equation*}
where $N\in\N^+$, $d\geq3$, $0<c\leq c_*:=\f{(d-2)^2}{4}$,  $c_*>0$ is the best constant of Hardy's inequality, and $V\in L^{s+\f d2}(\R^d)$ holds for $s>0$. Here $\lambda_i\big(-\Delta-{c}{|x|^{-2}}-V\big)$ denotes the $i$-th min-max level of Hardy-Schr\"odinger operator $H_{c,V}:=-\Delta-{c}{|x|^{-2}}-V $ in $\R^d$, which equals to the $i$-th negative eigenvalue (counted with multiplicity) of $H_{c,V}$ in $\R^d$ if it exists, and vanishes otherwise. We analyze the existence and analytical properties of the optimizers for the above inequality.
\end{abstract}


\section{Introduction}
In order to investigate the stability of quantum matter,  the following Lieb-Thirring type inequality and its generalizations of all kinds  were studied extensively over the past few decades (cf. \cite{EF06,Frank09,Frank,FLS08,LT75,LT76}):
\begin{equation}\label{compri}
	\operatorname{Tr}\Big(-\Delta-\f{c}{|x|^2}-V\Big)_-^s\leq L_{s,d,c}\int_{\R^d}V_+^{s+\f d2}dx,
\end{equation}
where $d\ge 3$, $c\in[0,c_*]$, $s\ge 0$ and $V\in L^{s+\f{d}{2}}(\R^d)$ is a real-valued potential. Here $V_+:=\max\{0,V\}$ denotes the positive part of $V$, and $c_*:=\f{(d-2)^2}{4}>0$ is the best constant of the following Hardy's inequality
\begin{equation}\label{hdi}
	\int_{\R^d}\f{c_*}{|x|^2}|u|^2dx\leq\int_{\R^d}|\n u|^2dx,\ \,\forall \,u\in  H^1(\R^d, \mathbb{C}).
\end{equation}
The inequality \eqref{compri} gives an upper bound,  in terms of the $L^{s+\f d2}$-norm of the potential $V_+$ in $\R^d$, on the magnitude  of the negative eigenvalues for the following Hardy-Schr\"odinger operator
$$H_{c, V}:=-\Delta-c|x|^{-2}-V\ \  \mbox{in}\,\ \R^d.$$
By the dual principle (cf. \cite{GLN21}), see also \eqref{dual} below, the existence and uniqueness of the first eigenfunction for the operator $H_{c, V}$ in $\R^d$ were well studied in \cite{KMK,LT76,Nam21,TGG} for any $c\in[0,c_*]$. In spit of this fact, as far as we know, there exist very few analytical results of other eigenfunctions for the operator $H_{c, V}$ in $\R^d$.

It deserves to point out that the inequality \eqref{compri} with $c=0$ is the well-known Lieb-Thirring   inequality, which was proposed originally by  Lieb and  Thirring in 1975 (cf. \cite{LT75,LT76}). The Lieb-Thirring  inequality plays a crucial role in mathematical physics, see \cite{LS10,LT75,LT76} and the references therein. Particularly, it follows from \cite{Frank,LT75,LT76} that the range of $s\geq0$ for the Lieb-Thirring  inequality satisfies
\begin{equation}\label{frlti:M}
	s\begin{cases}
		\geq1/2,\ \, &\text{if } \ d=1;\\ >0,\quad &\text{if }\ d=2;\\ \geq0, \quad &\text{if }\ d\geq3.
	\end{cases}
\end{equation}
However, there still exist (cf. \cite{Frank})  many challenging unsolved problems for the Lieb-Thirring  inequality. For this reason, the authors in \cite{FGL21,FGL24,GLN21} suggested a new scenario of analyzing the Lieb-Thirring inequality, for which the following finite-rank Lieb-Thirring inequality was studied:
\begin{equation}\label{frlti}
	\sum_{i=1}^N\left|\lambda_i\left(-\Delta-V\right)\right|^s\leq L_{s,d}^{(N)}\int_{\R^d}V_+^{s+\f d2}dx.
\end{equation}
Here $\lambda_i\left(-\Delta-V\right)$ denotes the $i$-th min-max level of the operator $-\Delta-V$ in $\R^d$, which equals  (cf. \cite[Theorem 11.6]{Lieb01}) to the $i$-th negative eigenvalue (counted with multiplicity) of the operator $-\Delta-V$ in $\R^d$ if it exists, and vanishes otherwise.
The existence of optimizing potentials $V$ for \eqref{frlti} was proved in \cite{FGL21,FGL24} for each $N\in \mathbb{N}^+$. Moreover, it was shown in \cite{FGL24} that if $s >\max\{0,2-d/2\}$, then the optimal Lieb-Thirring constant $L_{s,d}^{(N)}$ cannot be stationary for all $N\in \mathbb{N}^+$.

When $c=c_*$, Ekholm and Frank proved  in \cite{EF06} that the inequality \eqref{compri} holds for $d\geq3$ and $s>0$, which is completely different from the case $c=0$, see (\ref{frlti:M}) and \cite{Frank24} for more details.
More precisely, the authors in \cite{EF06} addressed the following inequality
\begin{equation}\label{chlt}
	\operatorname{Tr}\Big(-\Delta-\f{c_*}{|x|^2}-V\Big)_-^s\leq C_{s,d}\int_{\R^d}V_+^{s+\f d2}dx,\quad d\geq3\text{ \,and\, } s>0,
\end{equation}
where $C_{s,d}>0$ is the best constant.
It was shown in  \cite{EF06} that when $c=0$ and  $d\geq 3$,  a sufficiently weak potential $V$ may not give any bound state of \eqref{compri}. For example, if we take $V=\f{c_*}{|x|^2}$, then the inequality \eqref{hdi} immediately yields that the left-hand side of \eqref{compri} is zero. However, the right-hand side of \eqref{compri} is infinite, which further implies that the inequality \eqref{compri} with $c=0$ does not provide a meaningful bound for sufficiently weak potentials $V$. On the other hand, it was noted from  \cite{Ti99} that for $V\gneqq 0$, the operator $-\Delta-V-\f{c_*}{|x|^2}$ in $\R^d$ ($d\geq3$)  has always a bound state, but the integral on the  right-hand side of \eqref{chlt} can be made arbitrarily small. For this reason, Ekholm and Frank introduced in \cite{EF06} the concept of the  virtual level, where the analysis on the sharp value of the best constant $C_{s,d}$ for \eqref{chlt} was also imposed as a challenging open problem.

Stimulated by above facts, in this paper we are interested in the following finite-rank Hardy-Lieb-Thirring inequality
\begin{equation}\label{hlt}
	\sum_{i=1}^N\Big|\lambda_i\Big(-\Delta-\f{c_*}{|x|^2}-V\Big)\Big|^s\leq C_{s,d}^{(N)}\int_{\R^d}V_+^{s+\f d2}dx,
\end{equation}
where $N\in\N^+$,  $s>0$ and $c_*:=\f{(d-2)^2}{4}>0$ holds for $d\geq 3$. Here the $i$-th min-max level is defined by
\begin{equation*}
	\begin{split}
		&\lambda_i\Big(-\Delta-\f{c_*}{|x|^2}-V\Big)\\
:=&\inf_{\phi_1,\cdots,\phi_i\in\mathcal Q}\max\l\{\int_{\R^d}\Big(|\n \phi|^2-\f{c_*}{|x|^2}|\phi|^2-V|\phi|^2\Big)dx:\ \phi\in\text{Span}(\phi_1,\cdots,\phi_i),\,\|\phi\|_2=1\r\},
	\end{split}
\end{equation*}
where  $\phi_1,\cdots,\phi_i\text{ are orthonormal in } L^2(\R^d)$, and the Hilbert space $\mathcal Q$ is defined by
\begin{equation}\label{Q}
\mathcal Q:=\Big\{u\in L^2(\R^d, \mathbb{C}):\,\int_{\R^d}\Big(|\nabla u|^2-\f{c_*}{|x|^2}|u|^2\Big )dx<+\infty\Big\}
\end{equation}
equipped with the inner product
\begin{equation}\label{Qip}
	\big<u,v\big>_{ \mathcal Q}=\int_{\R^d}{\Big(\n \bar u\cdot\n  v-\f{c_*}{|x|^2}\bar u v+\bar u v\Big)dx}.
\end{equation}
One can check that the best constant $C_{s,d}^{(N)}>0$ of (\ref{hlt}) is non-decreasing with respect to $N\in\N^+$, and
\begin{equation}\label{CCN}
	C_{s,d}=\sup_{N\in\N^+}C_{s,d}^{(N)}=\lim_{N\to\infty}C_{s,d}^{(N)}<\infty,
\end{equation}
where $C_{s,d}>0$ is the best constant of \eqref{chlt}. Moreover,  if the potential $V$ is replaced by $V_+$, then the left-hand side of \eqref{hlt} increases. Therefore, without loss of generality, we may assume that $V(x)\geq0$.

We also have
\begin{equation*}
	t^{2s+d}\int_{\R^d} V(tx)^{s+\f d2}dx=t^{2s}\int_{\R^d}V(x)^{s+\f d2} dx,
\end{equation*}
and

\begin{equation*}
	\begin{split}
		&\inf_{\substack{u\in \mathcal{Q},\\ \|u\|_2=1}}\int_{\R^d}\Big( |\n u|^2-\f{c_*}{|x|^2}|u|^2-t^2V(tx)|u|^2\Big)dx\\
		=&t^2\inf_{\substack{u\in \mathcal{Q},\\ \|u\|_2=1}}\int_{\R^d}\Big(|\n u|^2-\f{c_*}{|x|^2}|u|^2-V(x)|u|^2\Big)dx,
	\end{split}
\end{equation*}
which imply that the inequality \eqref{hlt} is invariant under the scaling $t^2V(t\cdot)$  with $t>0$. Here $\mathcal Q $ is defined by \eqref{Q} below. Thus, one may assume that $V$ is normalized in $L^{s+\f d2}{(\R^d)}$, so that the inequality \eqref{hlt} can be rewritten equivalently as the following form:
\begin{equation}\label{sup}
	C_{s,d}^{(N)}=\sup_{\substack{0\leq V\in L^{s+\f d2}(\R^d),\\
			\int_{\R^d}V^{s+\f d2}dx=1}}\sum_{i=1}^N\Big|\lambda_i\Big(-\Delta-\f{c_*}{|x|^2}-V\Big)\Big|^s,\ N\in\N^+,
\end{equation}
where $s>0$ and $c_*=\f{(d-2)^2}{4}>0$ holds for $d\geq 3$. Therefore, in this paper we mainly analyze the maximization problem \eqref{sup}, which then presents equivalently the analysis of the inequality \eqref{hlt}. As discussed in Remark 1.1 below, all main results of the present paper still hold true, if the critical constant $c_*$ of \eqref{hlt} and (\ref{sup}) is replaced by any $c\in [0, c_*]$.

We now denote $0\leq\gamma=\gamma^*$ to be a non-negative self-adjoint  operator satisfying $\operatorname{Rank}(\gamma)\leq N,$ which has the form $\gamma=\sum_{i=1}^Nn_i|u_i \rangle \langle u_i|$ in terms of an orthonormal family $\{u_1,\cdots,u_N\}$  in $L^2(\R^d, \mathbb{C})$. Following \cite{Sim05}, the $q$-th Schatten norm of the operator $\gamma$ is defined by
\begin{equation*}
	\|\gamma\|_q:=\big(\operatorname{Tr}|\gamma|^q\big)^{1/q}=\big(\sum_{i=1}^Nn_i^q\big)^{1/q},\,\ 1\leq q<\infty,
\end{equation*}
and $\|\gamma\|_\infty:=\|\gamma\|$ denotes  the norm of the operator $\g$. The density function $\rho_\gamma$ of $\gamma$ is then defined by
\begin{equation*}
	\rho_\gamma:=\sum_{i=1}^Nn_i|u_i|^2 \in L^1(\R^d, \mathbb{R}).
\end{equation*}
We then set
\begin{equation*}
	\operatorname{Tr}\l(-\Delta-\f{c_*}{|x|^2}\r)\gamma:=\sum_{i=1}^Nn_i\int_{\R^d}\Big(|\nabla u_i|^2-\f{c_*}{|x|^2}|u_i|^2\Big)d x,
\end{equation*}
where   $\operatorname{Tr}\l(-\Delta-\f{c_*}{|x|^2}\r)\gamma=+\infty$ is assumed, if $u_i\notin \mathcal Q$ holds for some $i\in\N^+$. Here $\mathcal Q$ is defined by \eqref{Q}.

For $d\geq 3$ and $1< p\leq 1+ 2/d$, define
\begin{equation*}
	q:=\left\{
	\begin{aligned}
	&\f{2p+d-dp}{2+d-dp},\,\ \mbox{if}\,\ 1< p<1+\f 2d;\\
	&+\infty,\ \qquad \quad\mbox{if}  \ \ p=1+\f 2d.
	\end{aligned}\right.
\end{equation*}
We denote  $D_{p,d}^{(N)}$ the sharp constant of the following inequality
\begin{equation}\label{dual} D_{p,d}^{(N)}\|\rho_\gamma\|_{L^p(\R^d)}^{\f{2p}{d(p-1)}}\leq\|\gamma\|_q^{\f{p(2-d)+d}{d(p-1)}}\operatorname{Tr}\l(-\Delta-\f{c_*}{|x|^2}\r)\gamma,
\end{equation}
where as before $0\leq\gamma=\gamma^*$ is a non-negative self-adjoint  operator satisfying $\operatorname{Rank}(\gamma)\leq N$. Similar to \cite[Lemma 5]{FGL21}, one can show that the best constant $C_{s,d}^{(N)}>0$ of  \eqref{hlt} satisfies
\begin{equation}\label{CD}
	D_{p,d}^{(N)}\big(C_{s,d}^{(N)}\big)^{\f 2d}=\left(\f{2s}{2s+d}\right)^\f{2s}{d} \f{d}{2s+d} ,	
\end{equation}
where $1\leq N<\infty$, $d\geq3$, $p=\f{2s+d}{2s+d-2}$, and
\begin{equation}\label{duall}
s=\left\{
\begin{aligned}
\geq 1,\quad N\geq2;\\
>0,\quad N=1.
\end{aligned}\right.
\end{equation}
Applying \cite[Theorem 1.1]{EF06}, we then obtain from \eqref{CD} that $0<D_{p,d}^{(N)}<\infty$ holds for all $N\in\N^+$.

Especially, when $N=1$, the term $\|\gamma\|_q^{\f{p(2-d)+d}{d(p-1)}}$ in \eqref{dual} can be  eliminated, and the dual argument then implies from \eqref{dual} and \eqref{CD} that
\begin{equation}
C_{s,d}^{(1)}=C_{HGN}^{-\f{2m}{m-2}}\left(\f{2s}{2s+d}\right)^s\left(\f{d}{2s+d}\right)^{\f d2},\quad 2<m=\f{2d+4s}{d-2+2s}<2^*=\f{2d}{d-2},
\end{equation}
where $d\geq3$, $s>0$, and $C_{HGN}>0$ is the sharp constant of the inequality
\begin{equation}\label{gnh}
C_{HGN}^{\f{4m}{d(m-2)}}\Big(\int_{\R^d}|u|^{m}dx\Big)^{\f4{(m-2)d}}\leq\Big(\int_{\R^d}|\nabla u|^2-\f{c_*}{|x|^2}|u|^2dx\Big)\Big(\int_{\R^d}|u|^2\Big)^{\f{2m-md+2d}{(m-2)d}},\  u\in    \mathcal Q.
\end{equation}
Moreover, using a similar argument as in \cite{LT76}, one can deduce that for $d\geq3$ and $s>0$, $C_{s,d}^{(1)}$ always admits   optimizers of the form
\begin{equation*}
V(x)=t^2Q(tx)^{m-2}, \ \ t>0,	
\end{equation*}
where $Q\in\mathcal Q$ is (cf.  \cite[Theorem 1]{Nam21}) the unique radial positive solution of
\begin{equation}
		(-\Delta-c_*|x|^{-2})Q-Q^{m-1}+Q=0\text{ \ in \ }\R^d,
\end{equation}
see also \cite{Frank,LT76} for similar discussions on the Lieb-Thirring inequality \eqref{frlti} with $N=1$. Generally, it was shown in \cite{Nam21} that for any $N\in\N^+$, the quadratic form domain of the operator $-\Delta-c_*|x|^{-2}$ in $\R^d$ is $ \mathcal Q$,  which satisfies
\begin{equation*}
H^1(\R^d,\C)\subsetneq\mathcal Q\subset L^{2}(\R^d,\C).
\end{equation*}
We also get from \cite[Theorem 1.2]{Frank09} that $   \mathcal Q$ can be continuously embedded into  the fractional Sobolev space, $i.e.,$
\begin{equation}\label{fh}
\mathcal{Q}\subset H^t(\R^d,\C),\quad\forall \,t\in(0,1),
	\end{equation}
and hence
\begin{equation}\label{ce}    \mathcal{Q} \text{ can be compactly embedded into \ }L^{q}_{loc}(\R^d,\C) \text{ \ for all \ }2\leq {q}<2^*,
\end{equation}
see \cite[Theorem 1.4]{BRS} for more details.

The first result of this paper is concerned with the following existence of optimizers for \eqref{sup}, where $N\in\N^+$ is arbitrary.

\begin{Theorem}\label{exist}
For any fixed  $d\geq 3$, $N\in\N^+$ and $s>0$, let the problem $C_{s,d}^{(N)}$ be defined by \eqref{sup}. Then the problem $C_{s,d}^{(N)}$ admits at least one optimizer.
 \end{Theorem}
 Even though the similar existence of Theorem \ref{exist}  was proved earlier in some existing works, see \cite{FGL21,FGL24} and the references therein, the presence of the critical Hardy  term $c_*|x|^{-2}$ leads to some extra difficulties in the proof of Theorem \ref{exist}. Firstly, the critical Hardy  term $c_*|x|^{-2}$ yields the lack of translation invariance for the problem  $C_{s,d}^{(N)}$, and hence the existence method of \cite{FGL24} is not applicable to the problem  $C_{s,d}^{(N)}$. Secondly, once the critical Hardy  term $c_*|x|^{-2}$ is involved, then the eigenfunctions of the operator $-\Delta-{c_*}{|x|^{-2}}-V$ in $\R^d$ do not belong to $H^1(\R^d,\C)$. To overcome these difficulties, inspired by \cite[Theorem 13]{Nam21},   we first observe from \cite[Corollary 2]{FGL24} that the finite-rank Lieb-Thirring constant $L_{s,d}^{(N)}$ of \eqref{frlti} is   attainable, which helps us establish successfully in Section \ref{S2} the tightness of the maximizing sequence $\{V_n\}$ for $C_{s,d}^{(N)}$, in view of the analytical result $C_{s,d}^{(N)}>L_{s,d}^{(N)}$, see \eqref{gapcl}. Different from \cite{Nam21}, we however need to obtain the dichotomy of the eigenvalues in our case, see Lemma \ref{ed} for more details. We further  prove in Section \ref{S3} that any optimizing sequence of the problem \eqref{sup} possesses a subsequence that converges to an optimizer of $C_{s,d}^{(N)}$ strongly in $L^{s+\f d2}(\R^d,\R)$.

In 2006, Ekholm and Frank imposed  in \cite{EF06} the sharp value of the optimal constant $C_{s,d}$ defined by \eqref{CCN} as a challenging open problem. As far as we know, there is however no any progress on this open problem. Since the finite-rank Hardy-Lieb-Thirring inequality (\ref{hlt}) does not admit the translation invariance, the strict monotonicity of $C_{s,d}^{(N)}$ in $N\in\N^+$ is still unknown. Consequently, one may wonder whether the sharp constant $C_{s,d}$ can be attained by a finite-rank potential, $i.e.$, whether there exists a positive integer $N$ such that $C_{s,d}=C_{s,d}^{(N)}$. Therefore, Theorem \ref{exist} can help us understand the above open problem of \cite{EF06}.

 For any $N\in\N^+$ and $d\geq 3$, suppose  $V_N$ is an optimizer of $C_{s,d}^{(N)}$ defined by \eqref{sup}. Denote
 \begin{equation}\label{1:MM}
 	\lambda_i:=\lambda_i\big(H_{c_*, V_N}\big),\ \, H_{c_*, V_N}=-\Delta-V_N-{c_*}{|x|^{-2}}
 \end{equation}
 the corresponding $i$-th min-max level of the operator $H_{c_*, V_N}$ in $\R^d$, which equals  (cf. \cite[Theorem 11.6]{Lieb01}) to the $i$-th negative eigenvalue (counted with multiplicity) of the operator $H_{c_*, V_N}$ in $\R^d$ if it exists, and vanishes otherwise. Let $M\in\N^+$ be the number of negative eigenvalues (counted with multiplicity) for the operator $H_{c_*, V_N}$ in $\R^d$, $i.e.,$ $M\in\N^+$ is the smallest integer so that
 \begin{equation}\label{1:M}
 	\lambda_{M+1}=0 \ \, \mbox{and}\,\ \lambda_{i}<0 \ \, \mbox{holds for all}\,\ 1\le i\le M.
 \end{equation}
In this paper, we shall address the following analytical properties of optimizers for \eqref{sup}, where $N\in\N^+$ is arbitrary.

\begin{Theorem}\label{analytic}
For any fixed  $d\geq 3$, $N\in\N^+$ and $s>0$, let the problem $C_{s,d}^{(N)}$ be defined by \eqref{sup}, and let $V_N$ be an optimizer of $C_{s,d}^{(N)}$. Suppose   $M\in\N^+$ is the number of negative eigenvalues (counted with multiplicity) of the operator $H_{c_*, V_N}$ defined by (\ref{1:MM}) in $\R^d$. Then there exists a sufficiently large $R>0$ such that
\begin{equation}\label{1.20}
0<V_N(x)\leq C\l(\e ^{-\f{\sqrt{|\lambda_{N'}|}}{2}|x|}\r)^{\f 4{2s+d-2}}\text{ \quad for\quad }|x|>R,\ N'=\min\{M,N\},
\end{equation}
and
\begin{equation}\label{ev}
V_N=\l(\f{2s}{(d+2s)C_{s,d}^{(N)}}\sum_{i=1}^{N'}|\lambda_{i}|^{s-1}|u_i|^2\r)^{\f{2}{2s+d-2}}\in C^2(\R^d\setminus\{0\}),\,\ N'=\min(M,N),
\end{equation}
where $u=(u_1,\cdots, u_{N'})$ satisfying $u_i\in C^2(\R^d\setminus\{0\})$ solves  the following orthonormal system:
\begin{equation} \label{A:Gap} \Big[-\Delta-\f{c_*}{|x|^2}-\l(\f{2s}{(d+2s)C_{s,d}^{(N)}}\sum_{i=1}^{N'}|\lambda_{i}|^{s-1}|u_i|^2
\r)^{\f{2}{2s+d-2}}\Big]u_i=\lambda_iu_i \text{ \ in\ }\, \R^d,\,\  i=1,\cdots N'.
\end{equation}
Moreover, if $\lambda_{N+1}<0$ (namely, $M>N$), then $\lambda_N$ satisfies
\begin{equation}\label{Gap}
\lambda_N<\lambda_{N+1}<0.
\end{equation}
\end{Theorem}

The proof of Theorem \ref{analytic}  relies on the perturbation theory of linear operators (cf. \cite{Kat95}). Specially, analyzing the gap \eqref{Gap} is crucial for the proof of \eqref{ev}. Hence, we shall first prove \eqref{Gap} by precisely analyzing the eigenvalues of a perturbed operator, see Lemma \ref{L3.3} for more details.
We also note from the inequality \eqref{Gap} that there are ``no unfilled shells for the optimizers of the finite-rank Hardy-Lieb-Thirring inequality", in the spirit of what is known (cf. \cite{Lieb94}) for Hartree-Fock theory.

\begin{Remark}
Theorems \ref{exist} and \ref{analytic} remain valid if the critical constant $c_*$ of \eqref{sup} is replaced by any constant $c\in(0,c_*]$. In fact, one can check that if the critical constant $c_*$ of \eqref{sup} is replaced by any constant $c\in(0,c_*)$, then  the corresponding eigenfunctions belong to $H^1(\R^d)$, instead of $\mathcal Q$ defined in (\ref{Q}). This further implies that   Theorems \ref{exist} and \ref{analytic} still hold by the simpler analysis.
\end{Remark}

This paper is organized as follows. In Section \ref{S2}, we shall prove Proposition \ref{T2.1} on the tightness of the optimizing sequence for $C_{s,d}^{(N)}$ defined by \eqref{sup}. In Section \ref{S3}, we establish  Theorem \ref{exist}  on the existence of optimizers for $C_{s,d}^{(N)}$. Theorem \ref{analytic}  is then proved in Section \ref{S4}, which is concerned with the analytical properties of optimizers for  $C_{s,d}^{(N)}$.

\section{Tightness of the Optimizing Sequence}\label{S2}

For any fixed $s>0$, $d\geq3$ and $N\in\N^+$, let the problem  $C_{s,d}^{(N)}$ be defined by \eqref{sup}. From now on,  we always denote $\lambda_i(\mathcal L_n)$ to be the $i$-th min-max level of the operator $\mathcal L_n$ in $\R^d$. The main purpose of this section is to prove the following tightness of the optimizing sequence for $C_{s,d}^{(N)}$.

\begin{Proposition}\label{T2.1}
Let $\{V_n\geq0\}\subset L^{s+\f d2}(\R^d, \R)$ be a normalized optimizing sequence of $C_{s,d}^{(N)}$ defined by \eqref{sup}, where $s>0$, $d\geq3$ and $N\in\N^+$. Then the sequence $\{V_n\}$ is tight in the sense that
\begin{equation}\label{am}
 \lim_{R\to\infty}\liminf_{n\to\infty}\int_{|x|\leq R}V_n^{s+\f d2}dx= \lim\limits_{n\to\infty}\int_{\R^d}V_n^{s+\f d2}dx.
\end{equation}
\end{Proposition}

We remark that Proposition \ref{T2.1} is crucial for establishing Theorem \ref{exist}  on the existence of optimizers for  $C_{s,d}^{(N)}$. To establish Proposition \ref{T2.1}, it is necessary to investigate the impact of the maximizing sequence on the limiting behavior of the eigenvalues. Inspired by \cite[Lemma 17]{FGL24}, we first establish the following dichotomy of the spectrum.

\begin{Lemma}\label{ed}
For any fixed $s>0$ and $d\geq 3$, suppose the sequence $\{V_n\geq0
\}$ is  bounded uniformly in $L^{s+\f d2}(\R^d,\R)$. If there exists a sequence $\{R_n\}\subset\R$ satisfying $R_n\to\infty$ as $n\to\infty$ such that
\begin{equation}\label{0}
\lim_{n\to\infty}\int_{B_{2R_n}\setminus B_{R_n}}V_n^{s+\f d2}dx=0,
\end{equation}
then for any $N\in\N^+$, there exists $M\in\{0,\cdots,N\}$ such that up to a subsequence if necessary,
\begin{align}\label{2.3}
\sum_{i=1}^N\Big|\lambda_i\Big(-\Delta-\f{c_*}{|x|^2}-V_n\Big)\Big|^s=&\sum_{i=1}^M\Big|\lambda_i\Big(-\Delta-V_n\mathds{1}_{B_{R_n}}-\f{c_*}{|x|^2}\Big)\Big|^s\\
&+\sum_{i=1}^{N-M}\Big|\lambda_i\Big(-\Delta-V_n\mathds{1}_{\R^d\setminus B_{2R_n}}-\f{c_*}{|x|^2}\Big)\Big|^s+o(1)\text{ as }n\to\infty.\nonumber
	\end{align}
\end{Lemma}

\noindent{\bf Proof.}  Define
\begin{equation*}
A:=\left(\begin{matrix}
-\Delta-V_n\mathds{1}_{B_{R_n}}-\f{c_*}{|x|^2}&0\\
0&-\Delta-V_n\mathds{1}_{\R^d\setminus B_{2R_n}}-\f{c_*}{|x|^2}
\end{matrix}\right)
\end{equation*}
as an operator acting on $L^2(\R^d,\mathbb{C})\times L^2(\R^d,\mathbb{C})$, where the sequence $\{R_n\}$ satisfying $R_n\to\infty$ as $n\to\infty$ is given by \eqref{0}. One can   verify that
 \begin{equation*}
	\sigma_p^-(A)=\sigma_p^-\big(	-\Delta-V_n\mathds{1}_{B_{R_n}}-{c_*}{|x|^{-2}}\big) \cup\sigma_p^-\big(-\Delta-V_n\mathds{1}_{\R^d\setminus B_{2R_n}}-{c_*}{|x|^{-2}}\big),
\end{equation*}
where  $\sigma_p^-(\mathcal L)$ denotes  the set of all negative eigenvalues for the operator $\mathcal L$ in $\R^d$.
Thus,  in order to complete the proof of Lemma \ref{ed}, it suffices to prove the following stronger conclusion:
for any $i\in\N^+$,
\begin{equation*}\label{spectrum}
	\lambda_i\left(-\Delta-\f{c_*}{|x|^2}-V_n\right)=\lambda_i(A)+o(1)\ \ \ \text{as}\ \ n\to\infty.
\end{equation*}
Similar to \cite[Lemma 17]{FGL24}, one can prove that  for every $i\in\N^+$,
\begin{equation*}\label{spectrum1}
\lambda_i\left(-\Delta-\f{c_*}{|x|^2}-V_n\right)\leq\lambda_i(A)+o(1)\ \  \text{ as }\  n\to\infty.
\end{equation*}
It hence remains to prove the following reverse inequality
\begin{equation}\label{spectrum2}
	\lambda_i\left(-\Delta-\f{c_*}{|x|^2}-V_n\right)\geq\lambda_i(A)+o(1)\ \  \text{ as }\ n\to\infty,\ \forall i\in\N^+.
\end{equation}

	 Let $f\in C_c^\infty(\R^d,[0,1])$ be a radial function such that $f(x)\equiv1$ on $B_1$ and $f(x)\equiv0$ on $B_{5/4}^c$, and let $g\in C^\infty(\R^d,[0,1])$ be a radial function such that $g(x)\equiv1$ on $B_2^c$ and $g(x)\equiv0$ on $B_{7/4}$. For the sequence $\{R_n\}\subset\R$ satisfying \eqref{0}, define
	\begin{equation}\label{a.1}
		f_n(x):=f\big({x}/{R_n}\big),\ \ \ g_n(x):=g\big({x}/{R_n}\big), \ \ \  h_n(x):=h\big({x}/{R_n}\big),
	\end{equation}
where $h=\sqrt{1-f^2-g^2}$. Using the  IMS localization formula (cf. \cite[Theorem 3.2]{CFKS87}), we then obtain that for any $\epsilon_n\in(0,1)$,
	\begin{align}\label{c.1}
&-\Delta-\f{c_*}{|x|^2}-V_n\nonumber\\
=&(1-\epsilon_n)\left[
f_n\Big(-\Delta-V_n\mathds{1}_{B_{R_n}}-\f{c_*}{|x|^2}\Big)f_n+g_n\Big(-\Delta-V_n\mathds{1}_{B_{2R_n}^c}-\f{c_*}
{|x|^2}\Big)g_n\right. \\
&\qquad\qquad\,\left.+h_n\Big(-\Delta-\f{c_*}{|x|^2}\Big)h_n-|\n f_n|^2-|\n g_n|^2-|\n h_n|^2\right]\nonumber\\
&+\epsilon_n\Big(-\Delta-\f{c_*}{|x|^2}-\epsilon_n^{-1}V_n\mathds{1}_{B_{2R_n}\setminus B_{R_n}}-V_n\mathds{1}_{ B_{R_n}}-V_n\mathds{1}_{B_{2R_n}^c}\Big)\ \ \text{ in }\, \R^d.\nonumber
\end{align}
In view of \eqref{0}, choose $\epsilon_n:=\|V_n\mathds{1}_{B_{2R_n}\setminus B_{R_n}}\|_{s+\f d2}=o(1)$ as $n\to\infty$. By H\"older's inequality and the uniform boundedness of $\{V_n\}$ in $L^{s+\f d2}(\R^d,\R)$, one can calculate from \eqref{gnh} that for any  $u\in    \mathcal Q$ satisfying $\|u\|_{2}=1$,
\begin{equation*}
\begin{split}
&\quad\int_{\R^d}\Big[|\n u|^2-\f{c_*}{|x|^2}|u|^2-\epsilon_n^{-1}V_n\mathds{1}_{B_{2R_n}\setminus B_{R_n}}|u|^2-\Big(V_n\mathds{1}_{ B_{R_n}}+V_n\mathds{1}_{B_{2R_n}^c}\Big)|u|^2\Big]dx\\
&\geq C_1\|u\|_m^{\f{4m}{d(m-2)}}-\epsilon_n^{-1}\big\|V_n\mathds{1}_{B_{2R_n}\setminus B_{R_n}}\big\|_{s+\f d2}\, \|u\|_{m}^{2}-\|V_n\|_{s+\f d2}\, \|u\|_{m}^{2}\\
&=C_1\|u\|_m^{\f{4m}{d(m-2)}}-\|u\|_{m}^{2}-\|V_n\|_{s+\f d2}\|u\|_{m}^{2}\\
&\geq -C_2\ \ \text{ as }\, n\to\infty,
\end{split}
\end{equation*}
where $C_2>0$, $m=\f{2d+4s}{d-2+2s}>0$,  and the last inequality follows from the fact that $\f{4m}{d(m-2)}>2$. We hence obtain that for sufficiently large $n>0$,
\begin{equation}\label{c.2}
-\Delta-\f{c_*}{|x|^2}-\epsilon_n^{-1}V_n\mathds{1}_{B_{2R_n}\setminus B_{R_n}}-V_n\mathds{1}_{ B_{R_n}}-V_n\mathds{1}_{B_{2R_n}^c}\geq -C_2\ \ \text{ in }\, \R^d.
\end{equation}	
Additionally, it yields  from Hardy's inequality \eqref{hdi} that
	\begin{equation}\label{c.3}
 h_n\Big(-\Delta-\f{c_*}{|x|^2}\Big)h_n\geq 0\ \ \text{ in }\, \R^d.
		\end{equation}
We thus deduce from \eqref{c.1}--\eqref{c.3} that for sufficiently large $n>0$,
\begin{equation}\label{lb}
\begin{split}
&\, -\Delta-\f{c_*}{|x|^2}-V_n\\
\geq&\,
(1-\epsilon_n)\left[f_n\Big(-\Delta-V_n\mathds{1}_{B_{R_n}}-\f{c_*}{|x|^2}\Big)f_n+g_n\Big(-\Delta-V_n\mathds{1}_{B_{2R_n}^c}-\f{c_*}{|x|^2}\Big)g_n\right]\\
&-C(R_n^{-2}+\epsilon_n)\ \ \text{ in }\, \R^d,
\end{split}
\end{equation}
where $\f{1}{R_n}\to0$ and $\epsilon_n\to0$ as $n\to\infty$.

Define the partial isometry
\begin{equation*}
		\mathfrak{I}_n:\phi\in L^2(\R^d)\mapsto \left(\begin{matrix}
			f_n\phi\\g_n\phi\\h_n\phi
\end{matrix}\right)\in \big(L^2(\R^d)\big)^3.
\end{equation*}
We then rewrite \eqref{lb} as
\begin{equation*}
		\begin{split}
			-\Delta-\f{c_*}{|x|^2}-V_n\geq&(1-\epsilon_n)\mathfrak{I}_n^{-1}\left(\begin{matrix}
				-\Delta-V_n\mathds{1}_{B_{R_n}}-\f{c_*}{|x|^2}&0&0\\
				0&-\Delta-V_n\mathds{1}_{B_{2R_n}^c}-\f{c_*}{|x|^2}&0\\
				0&0&0
			\end{matrix}\right)\mathfrak{I}_n\\
			& -C(R_n^{-2}+\epsilon_n)\ \ \text{ in }\, \R^d,
		\end{split}
\end{equation*}
where $C>0$ is independent of $n>0$. Applying the min-max principle (cf. \cite[Theorem 11.8]{Lieb01}), we further conclude that for sufficiently large $n>0$,
	\begin{equation*}
		\begin{split}
			0&\geq\lambda_i\left(-\Delta-\f{c_*}{|x|^2}-V_n\right)\\
			&\geq(1-\epsilon_n)\lambda_i\left(\begin{matrix}
				-\Delta-V_n\mathds{1}_{B_{R_n}}-\f{c_*}{|x|^2}&0&0\\
				0&-\Delta-V_n\mathds{1}_{B_{2R_n}^c}-\f{c_*}{|x|^2}&0\\
				0&0&0
			\end{matrix}\right)-C(R_n^{-2}+\epsilon_n)\\
			&\geq\lambda_i\left(\begin{matrix}
				-\Delta-V_n\mathds{1}_{B_{R_n}}-\f{c_*}{|x|^2}&0\\
				0&-\Delta-V_n\mathds{1}_{B_{2R_n}^c}-\f{c_*}{|x|^2}
			\end{matrix}\right)-C(R_n^{-2}+\epsilon_n),
		\end{split}
	\end{equation*}
where $\f{1}{R_n}\to0$ and $\epsilon_n\to0$ as $n\to\infty$. This proves the reverse inequality \eqref{spectrum2},  and the proof of Lemma \ref{ed} is therefore complete.\qed

\subsection{Proof of Proposition \ref{T2.1}}

By applying Lemma \ref{ed}, we are now ready to address the proof of Proposition \ref{T2.1}.

\vspace{.1cm}
%Recall from \cite[Corollary 2]{FGL24} that the best constant $L_{s,d}^{(N)}$ of \eqref{frlti} is attainable, where $s>0,$ $d\geq3$ and $N\in\N^+$. Let $0\leq V\in L^{s+\f 2d}(\R^d)$ be an optimizer of $L_{s,d}^{(N)}$. It then follows from \eqref{sup1} that
%\begin{equation}\label{gapcl}
%	\begin{split}
%	C_{s,d}^{(N)}\geq \f{\sum_{i=1}^N\left|\lambda_i\left(-\Delta-\f{c_*}{|x|^2}-V\right)\right|^s}{ \int_{\R^d}V^{s+\f d2}dx}
%	>\f{\sum_{i=1}^N\big|\lambda_i\left(-\Delta-V\right)\big|^s}{ \int_{\R^d}V^{s+\f d2}dx}=L_{s,d}^{(N)}.
%	\end{split}
%\end{equation}

\noindent{\bf Proof of Proposition \ref{T2.1}.}
Let $\{V_n\geq0\}\subset L^{s+\f d2}(\R^d,\R)$ be an  optimizing sequence of $C_{s,d}^{(N)}$ defined by \eqref{sup}, where $s>0$, $d\geq3$ and $N\in\N^+$ are arbitrary. 	We shall carry out the proof by three steps.

\vspace{.1cm}
{\it Step 1.} We claim that
\begin{equation}\label{gig}
 \Lambda_i^D\leq0,\ \  \forall\ i\in\N^+,
\end{equation}
where the constant $\Lambda_i^D$ is defined by
\begin{equation}\label{gid}
	\lim_{n\to\infty} \lambda_i\Big(-\Delta-V_n{\mathds{1}_{ B_{2R_n}^c}}-\f{c_*}{|x|^2}\Big)
	_{B_{{7R_n}/{4}}^c}:=\Lambda_i^D,\ \  \forall\ i\in\N^+.
\end{equation}
Here $\{R_n\}\subset\R$ is any sequence satisfying $\lim\limits_{n\to\infty}R_n=\infty$, and the notation $\Big(-\Delta-V_n{\mathds{1}_{ B_{2R_n}^c}}-\f{c_*}{|x|^2}\Big)
_{B_{{7R_n}/{4}}^c}$ denotes the operator restricted  to the domain $B^c_{7R_n/4}$ with Dirichlet boundary condition.

Let $0<\mu_1<\mu_2\leq\cdots\leq\mu_i$ be the $i$-first eigenvalues of the operator $-\Delta_{B_2\setminus B_{7/4}}$ with Dirichlet boundary condition (cf. \cite[Section 6.5, Theorem 1]{Evans}), and suppose   $v_1,\cdots,v_i\in H^1_0(B_2\backslash B_{7/4},\, \R)$ are the associated orthonormal eigenfunctions, so that
\begin{equation}\label{213}
	\begin{cases}
		-\Delta v_j=\mu_j v_j&\text{ \ in }B_{2}\backslash B_{7/4},\\[1.5mm]
		v_j(x)=0  &\text{ \ on  }\partial(B_{2}\backslash B_{7/4}), \ \, j=1,\cdots,i.
	\end{cases}
\end{equation}
Set for $j=1,\cdots,i$,
\begin{equation}\label{214}
	v_{jn}:=\begin{cases}
	    R_n^{-\f{d}{2}}v_j\big(R^{-1}_nx\big)\ \text{ \ in\ }\, B_{2R_n}\setminus B_{7R_n/4},\\
    0 \qquad\qquad\qquad\,    \text{ \ in\ }\, B_{7R_n/4}\cup B^c_{2R_n},
	\end{cases}
\end{equation}
where  $R_n$ is as above.
Since $v_1,\cdots_,v_i$ form an orthonormal family in $L^2(\R^d,\R)$, we deduce  that
\begin{equation}\label{216}
	\langle v_{jn},v_{kn}\rangle=\delta_{jk}, \quad j,\,k=1,\cdots,i.
\end{equation}
By the min-max principle, one can calculate from \eqref{gid}--\eqref{216} that for any $i\in\N^+$,
\begin{equation}
\begin{split}
\Lambda_i^D&=\lim_{n\to\infty} \lambda_i\Big(-\Delta-V_n{\mathds{1}_{ B_{2R_n}^c}}-\f{c_*}{|x|^2}\Big)
_{B_{{7R_n}/{4}}^c}\\
&\leq \lim_{n\to\infty}\max_{\substack{\phi\in\text{span}(v_{1n},\cdots,v_{in})\\\|\phi\|_2=1}}
\int_{\R^d}\Big(|\n \phi|^2-V_n{\mathds{1}_{ B_{2R_n}^c}}|\phi|^2 -\f{c_*}{|x|^2}|\phi|^2\Big)dx\\
&\leq  \lim_{n\to\infty}\max_{\substack{\phi\in\text{span}(v_{1n},\cdots,v_{in})\\\|\phi\|_2=1}}\int_{\R^d}|\n \phi|^2 dx\\
		&\leq\lim_{n\to\infty}R_n^{-2}\sum_{j=1}^i\mu_j=0.
	\end{split}
\end{equation}
This thus yields that the claim \eqref{gig} holds true.

{\it Step 2.} In this step, we claim that for all $i\in\N^+$,
\begin{equation}\label{hardyv}
\lambda_i\Big(-\Delta-V_n{\mathds{1}_{ B_{2R_n}^c}}-\f{c_*}{|x|^2}\Big)
		_{B_{{7R_n}/{4}}^c}=\lambda_i\Big(-\Delta-V_n\mathds{1}_{ B_{2R_n}^c}-\f{c_*}{|x|^2}\Big)
		+o(1)\  \, \text{ as }\, n\to\infty,
\end{equation}
where $\{R_n\}\subset\R$ is any sequence satisfying $\lim\limits_{n\to\infty}R_n=\infty$.

Denote
\begin{equation}\label{gin}
	\lim_{n\to\infty}\lambda_i\Big(-\Delta-V_n\mathds{1}_{ B_{2R_n}^c}-\f{c_*}{|x|^2}\Big):=\Lambda_i,\ \ \ \forall\ i\in\N^+.
\end{equation}
The min-max principle then yields that
\begin{equation}\label{gig1}
	\Lambda_i^D\geq\Lambda_i,\ \ \ \forall\ i\in\N^+,
\end{equation}
where $\Lambda_i^D$ is defined in \eqref{gid}. If $\Lambda_{i^0}=0$ holds for some $i^0\in\N^+$, then it  follows immediately from  \eqref{gig} that  $\Lambda_i^D\leq\Lambda_i=0$ holds for  $i=i^0, i^0+1, \cdots$. Together with \eqref{gig1}, we thus conclude that \eqref{hardyv} holds true, and we are done. For this reason, without loss of generality, we next suppose that  $\Lambda_1\leq\cdots\leq\Lambda_N<0$, where $N\in\N^+$ is arbitrary.

%\textcolor{blue}{
%\begin{equation}\label{gid1}
%\lim_{n\to\infty}	\lambda_i\Big(-\Delta-V_n{\mathds{1}_{ B_{2R_n}^c}}-\f{c_*}{|x|^2}\Big)
	%	_{B_{{7R_n}/{4}}^c}:=\Lambda_i^D\leq 0,\ \ \ \forall\ i\
	%	\end{equation}
	%		and
			%\begin{equation}\label{gin1}
				%\lim_{n\to\infty}\lambda_i\Big(-\Delta-V_n\mathds{1}_{ B_{2R_n}^c}-\f{c_*}{|x|^2}\Big):=\Lambda_i\leq0,\ \ \ \forall\ i\in\N^+.
%				\end{equation}
%Note from the min-max principle that
%\begin{equation}\label{gig1}
%	\Lambda_i^D\geq\Lambda_i,\ \ \ \forall\ i\in\N^+.
%\end{equation}
%Thus, it remains to prove that
%\begin{equation}
%	\Lambda_i^D\leq\Lambda_i,\ \ \ \forall\ i\in\N^+.
%\end{equation}}

Since $\Lambda_1\leq\cdots\leq\Lambda_N<0$, we obtain from \eqref{gin} that  for sufficiently large $n>0$,
\begin{equation*}\lambda_{in}:=\lambda_i\left(-\Delta-V_n\mathds{1}_{ B_{2R_n}^c}-\f{c_*}{|x|^{2}}\right)<0,\ \ \ i=1, \cdots, N,
\end{equation*}
which yields that $\lambda_{in}$ is a negative eigenvalue, and thus there exists an eigenfunction associated to the eigenvalue $\lambda_{in}$.
%The same argument of \cite[Lemma 11.6]{Lieb01} yields that there exists an eigenfunction $u_{in}\in\mathcal {Q}$ corresponding to the eigenvalue $\lambda_{in}<0$.
Let $u_{in}$ be the $L^2$-normalized (real) eigenfunction associated to  $\lambda_{in}$, namely,
\begin{equation}\label{a.2}
\l(-\Delta-V_n\mathds{1}_{ B_{2R_n}^c}-\f{c_*}{|x|^2}\r)u_{in}=\lambda_{in}u_{in}\ \ \text{ in }\,\  \R^d,\ \ \ i=1, \cdots, N.
\end{equation}
Multiplying  \eqref{a.2} by $u_{in}$ and integrating over $\R^d$, it then follows  from \eqref{gnh} that
\begin{equation*}
\begin{split}
\int_{\R^d}\Big(|\n u_{in}|^2-\f{c_*}{|x|^2}| u_{in}|^2\Big)dx
&=\int_{\R^d}V_n\mathds{1}_{ B_{2R_n}^c}| u_{in}|^2dx+\lambda_{in}\\
&\leq \|V_n\mathds{1}_{ B_{2R_n}^c}\|_{s+\f d2}\|u_{in}\|_m^{2}-\lambda_{in}\\
&\leq C \left(\int_{\R^d}|\n u_{in}|^2-\f{c_*}{|x|^2}| u_{in}|^2dx\right)^{\f{d(m-2)}{2m}}+C,
\end{split}
\end{equation*}
where $2<m:=\f{2d+4s}{d-2+2s}<2^*=\frac{2d}{d-2}$, and $C>0$ is  independent of $n>0$. Since $0<{\f{d(m-2)}{2m}}<1$, we conclude from  \eqref{fh} that for $i=1, \cdots, N$,
\begin{align}\label{un}
\{u_{in}\}_n\ \text{ is bounded uniformly  in} \ \mathcal Q\cap L^q(\R^d,\R),\ \ \ 2\leq q<2^*.
\end{align}

Moreover, denote
\begin{equation}
\label{dln}l_n:=\sqrt{1-g_n^2},
\end{equation}
where $0\leq g_n\leq 1$ is as in \eqref{a.1}. It can be verified that $l_n\in C^\infty_c(\R^d, [0, 1])$ satisfies
\begin{equation}\label{2.14}
l_n(x)=1\ \ \text{for}\ \ |x|\leq 7R_n/4, \ \  l_n(x)=0\ \ \text{for}\ \ |x|\geq 2R_n.
\end{equation}
Multiplying  \eqref{a.2} by $l_n^2{u}_{in}$ and integrating over $\R^d$, we obtain  that for $i=1,\cdots,N$,
\begin{equation}\label{2.13}
\begin{split}
0&=\int_{\R^d}\Big(-\Delta u_{in}-V_n\mathds{1}_{ B_{2R_n}^c}u_{in}-\f{c_*}{|x|^2}u_{in}-\lambda_{in}u_{in}\Big)l_n^2{u}_{in}dx\\
&=\int_{\R^d}\n u_{in}\cdot\n(l_n^2{u}_{in})dx-\int_{\R^d}\Big(\f{c_*}{|x|^2}|l_nu_{in}|^2+\lambda_{in}|l_nu_{in}|^2\Big)dx\\
&=\int_{\R^d}\Big(|\n (l_nu_{in})|^2-\f{c_*}{|x|^2}|l_nu_{in}|^2\Big)dx+|\lambda_{in}|\int_{\R^d}|l_nu_{in}|^2dx-\int_{\R^d}|u_{in}|^2|\n l_n|^2dx.
\end{split}
	\end{equation}
We thus deduce  from \eqref{a.1}, \eqref{dln} and \eqref{2.13} that
\begin{equation*}
\begin{split}
0\leq&\int_{\R^d}\l(|\n (l_nu_{in})|^2-\f{c_*}{|x|^2}|l_nu_{in}|^2\r)dx+|\Lambda_i|\int_{\R^d}|l_nu_{in}|^2dx\\
=&\big(|\Lambda_i|-|\lambda_{in}|\big)\int_{\R^d}|l_nu_{in}|^2dx+\int_{\R^d}|u_{in}|^2|\n l_n|^2dx\\
\leq&|\Lambda_i-\lambda_{in}|+\f{C}{R_n^2}\to 0\ \ \ \text{as}\ \ n\to\infty,\ \ \ i=1, \cdots, N.
\end{split}
\end{equation*}
This further  yields that
	\begin{equation}\label{lu0}
		l_nu_{in}\to0\ \ \text{ strongly  in }    \mathcal{Q}\ \text{ as } \  n\to\infty,\ \ \ i=1, \cdots, N,
	\end{equation}
	where $\mathcal Q$ is defined by \eqref{Q}.

Define
\begin{equation}\label{dvn}v_{in}:=u_{in}-l_nu_{in},\ \ \ i=1, \cdots, N.
\end{equation}
One can verify from  \eqref{lu0} that
\begin{equation}\label{aor}
\lim\limits_{n\to\infty}\big<v_{in},v_{jn}\big>_2=\lim\limits_{n\to\infty}\big<u_{in},u_{jn}\big>_2=\delta_{ij},\ \ \ i, j=1, \cdots, N.
\end{equation}
Set
	\begin{equation}\label{c.6}
		B_n(u,v):=\Big\langle u,\ \Big(-\Delta-V_n{\mathds{1}_{ B_{2R_n}^c}}-\f{c_*}{|x|^2}\Big)v\Big \rangle_2.
	\end{equation}
We calculate  from \eqref{gin}, \eqref{a.2}, \eqref{2.14}  and \eqref{lu0} that for $i, j=1, \cdots, N$,
	\begin{equation}\label{c.4}
		\lim_{n\to\infty}B_n(u_{in},\,u_{jn})=\lim_{n\to\infty}\big<u_{in},\lambda_{jn}u_{jn}\big>_2=\delta_{ij}\Lambda_j,
	\end{equation}
	\begin{equation}
		\lim_{n\to\infty}B_n(l_nu_{in},\,l_nu_{jn})=\lim_{n\to\infty}\big(\big\langle l_nu_{in},\ l_nu_{jn}\big\rangle_{\mathcal Q}-\big\langle l_nu_{in},l_nu_{jn}\big\rangle_2\big)=0,
	\end{equation}
	and
	\begin{equation}\label{c.5}
		\begin{split}
			\l|B_n(u_{in},l_nu_{jn})\r|&=\big|\big<u_{in},l_nu_{jn}\big>_{\mathcal Q}-\big<u_{in},l_nu_{jn}\big>_2\big|\\[1mm]
			&\leq \big|\big<u_{in},l_nu_{jn}\big>_{\mathcal Q}\big|+\big|\big<u_{in},l_nu_{jn}\big>_2\big|\\[1mm]
			&\leq \|u_{in}\|_{\mathcal Q}\|l_nu_{jn}\|_{\mathcal Q}+\|u_{in}\|_2\|l_nu_{jn}\|_2\\[1mm]
			&\to0\text{\quad as }\, n \to\infty,
		\end{split}
\end{equation}
where we have used \eqref{un} and the Cauchy-Schwarz inequality. We thus deduce from \eqref{c.6}--\eqref{c.5} that for $ i, j=1, \cdots, N$,
\begin{equation}\label{c.06}
\begin{split}
B_n(v_{in},\,v_{jn})=&B_n(u_{in},\,u_{jn})+	B_n(l_nu_{in},\,l_nu_{jn})-2	B_n(v_{in},\,l_nu_{jn})\\
\to&\delta_{ij}\Lambda_j\text{\quad as}\,\
 n\to\infty,
\end{split}	
\end{equation}
which implies that
\begin{equation}\label{c.7}
	\begin{split}
		\Big\langle v_{in},\ \Big(-\Delta-V_n{\mathds{1}_{ B_{2R_n}^c}}-\f{c_*}{|x|^2}\Big)v_{jn}\Big \rangle_2\to\delta_{ij}\Lambda_j\text{\quad as }\ n\to\infty,\ \ i, j=1, \cdots, N.
	\end{split}
\end{equation}

Note from \eqref{2.14} and \eqref{dvn} that
\begin{align}\label{c.9}
	v_{in}(x)\equiv0\ \ \text{for} \ \ |x|\leq7R_n/4,\ \ i=1, \cdots,N.
\end{align}
 By  the min-max principle, we then deduce from \eqref{c.7} and \eqref{c.9} that
\begin{align}\label{2.31}
\Lambda_1=&\lim_{n\to\infty}\Big\langle v_{1n},\ \Big(-\Delta-V_n{\mathds{1}_{ B_{2R_n}^c}}-\f{c_*}{|x|^2}\Big)v_{1n}\Big \rangle_2\nonumber\\
=&\lim_{n\to\infty}\Big\langle v_{1n},\ \Big(-\Delta-V_n{\mathds{1}_{ B_{2R_n}^c}}-\f{c_*}{|x|^2}\Big)_{B_{7R_n/4}^c}v_{1n}\Big \rangle_2\\
\geq&\lim_{n\to\infty} \lambda_1\Big(-\Delta-V_n\mathds{1}_{ B_{2R_n}^c}-\f{c_*}{|x|^2}\Big)_{B_{7R_n/4}^c} \|v_{1n}\|_2^2=\Lambda_1^D,\nonumber
\end{align}
where the last identity follows from \eqref{gid} and \eqref{aor}.

Denote $ w_{1n}:=v_{1n}$. By  the Schmidt orthogonalization, we obtain from \eqref{c.9} that for every  $n\in\N^+$, there exists $c_{1n}\in\R$ such that the function 	
\begin{equation}\label{c.12}
	w_{2n}:=v_{2n}+c_{1n}w_{1n}\in\mathcal Q	\end{equation}
satisfies
\begin{equation}\label{c.11}
	w_{2n}(x)=0\ \ \text{for} \ \ |x|\leq7R_n/4,
\end{equation}
and
\begin{equation}\label{c.10}
	\langle w_{2n}, w_{1n}\rangle_2=0,
\end{equation}
where
\begin{equation}\label{c.10.5}
\lim_{n\to\infty}c_{1n}=0,
\end{equation}
due to the limit \eqref{aor}. We then get from \eqref{c.06} and \eqref{c.12} that	
\begin{equation}\label{2.33}
\begin{split}
			&\lim_{n\to\infty}B_n(w_{2n},w_{2n})\\
			=&\lim_{n\to\infty}\l[B_n(v_{2n},v_{2n})+2c_{1n}B_n(v_{1n},v_{2n})+c_{1n}^2B_n(v_{2n},v_{2n})\r]=\Lambda_2.
	\end{split}
\end{equation}
Together with \eqref{c.6} and \eqref{c.11}, it further gives that
\begin{equation}\label{2.35a}
\lim_{n\to\infty}\Big\langle w_{2n},\ \Big(-\Delta-V_n{\mathds{1}_{ B_{2R_n}^c}}-\f{c_*}{|x|^2}\Big)_{B_{7R_n/4}^c}w_{2n}\Big \rangle_2=\Lambda_2.
\end{equation}
Applying the min-max principle again, we conclude from \eqref{2.31}, \eqref{c.10} and \eqref{2.35a} that
\begin{equation}\label{2.41}
\sum_{i=1}^2\Lambda_i=	\lim_{n\to\infty}\sum_{i=1}^2\Big\langle w_{in},\ \Big(-\Delta-V_n{\mathds{1}_{ B_{2R_n}^c}}-\f{c_*}{|x|^2}\Big)_{B_{7R_n/4}^c}w_{in}\Big \rangle_i\geq\sum_{i=1}^2\Lambda_i^D.
\end{equation}
Since it follows from \eqref{gig1} and \eqref{2.31} that $\Lambda_1=\Lambda_1^D$, this further implies from \eqref{gig1} that
\begin{align*}
\Lambda_2=\Lambda_2^D.
\end{align*}

By induction, we thus obtain from above that
$$
\Lambda_i=\Lambda_i^D,\ \  \forall i=1,\cdots,N,
$$
$i.e., $ (\ref{hardyv}) holds true, and Step 2 is therefore proved.

\vspace{.1cm}

{\it Step 3.} In this step, we complete the proof of (\ref{am}).

We first claim that for any $i\in\N^+$,
\begin{equation}\label{0.10}
	\lambda_i\left(-\Delta-V_n\mathds{1}_{ B_{2R_n}^c}-\f{c_*}{|x|^2}\right)=\lambda_i\left(-\Delta-V_n{\mathds{1}_{ B_{2R_n}^c}}\right)+o(1)\ \ \text{ as }\, n\to\infty,
\end{equation}
where the sequence $\{R_n\}$ satisfies \eqref{0}.
Actually, one can check that for any $u\in\mathcal Q$ satisfying $u(x)=0$ for $|x|\leq{{7R_n}/{4}}$,
\begin{equation*}
\int_{\R^d}\f{c_*}{|x|^2}|u|^2dx=	\int_{B_{{7R_n}/{4}}^c}\f{c_*}{|x|^2}|u|^2dx\leq \f{C}{R_n^2}\int_{B_{{7R_n}/{4}}^c}|u|^2dx\to0\text{\quad as}\ \ n\to\infty,
	\end{equation*}
where the sequence $\{R_n\}$ satisfies \eqref{0}.
We then obtain from \eqref{hardyv} that for any $i\in\N^+$,
\begin{equation*}
\begin{split}
\lambda_i\Big(-\Delta-V_n\mathds{1}_{ B_{2R_n}^c}-\f{c_*}{|x|^2}\Big)&=\lambda_i\Big(-\Delta-V_n{\mathds{1}_{ B_{2R_n}^c}}-\f{c_*}{|x|^2}\Big)
			_{B_{{7R_n}/{4}}^c}
+o(1)\\
&=\lambda_i\left(-\Delta-V_n\mathds{1}_{ B_{2R_n}^c}\right)_{B_{{7R_n}/{4}}^c}
+o(1)\\
&=\lambda_i\left(-\Delta-V_n\mathds{1}_{ B_{2R_n}^c}\right)
+o(1)\ \ \text{ as }\, n\to\infty,
\end{split}
\end{equation*}
where the last identity follows from a similar approach  of \eqref{hardyv}. This proves the claim (\ref{0.10}).

We now define
\begin{equation}\label{2.43}
	\al:=	\lim_{R\to\infty}\liminf_{n\to\infty}\int_{|x|\leq R}V_n^{s+\f d2}dx\in[0,1].
\end{equation}
By \cite{Lions841}, there exists a sequence $\{R_n\}\subset\R$ satisfying $R_n\to\infty$ as $n\to\infty$ such that
\begin{equation}\label{a.4}
\lim_{n\to\infty}\int_{B_ {R_n}}V_n^{s+\f d2}dx=\al\ \ \text{ and }\	\lim_{n\to\infty}\int_{B_{2R_n}\setminus B_{R_n}}V_n^{s+\f d2}dx=0.
		\end{equation}
Recall from \cite[Corollary 2]{FGL24} that the best constant $L_{s,d}^{(N)}$ of \eqref{frlti} is attainable, where $s>0,$ $d\geq3$ and $N\in\N^+$. Let $0\leq V\in L^{s+\f d2}(\R^d)$ be an optimizer of $L_{s,d}^{(N)}$. It then follows from \eqref{hlt} that
\begin{equation}\label{gapcl}
\begin{split}
C_{s,d}^{(N)}\geq \f{\sum_{i=1}^N\left|\lambda_i\left(-\Delta-\f{c_*}{|x|^2}-V\right)\right|^s}{ \int_{\R^d}V^{s+\f d2}dx}
	>\f{\sum_{i=1}^N\big|\lambda_i\left(-\Delta-V\right)\big|^s}{ \int_{\R^d}V^{s+\f d2}dx}=L_{s,d}^{(N)}.
\end{split}
\end{equation}  	
Denote
\begin{equation}\label{e0}\epsilon_0:=C_{s,d}^{(N)}-L_{s,d}^{(N)}>0.
\end{equation}
By  Lemma \ref{ed}, we then deduce from \eqref{0.10} and \eqref{a.4} that for any $N\in\N^+$, there exists $M\in\{0,\cdots,N\}$ such that up to a subsequence if necessary,
\begin{align*}
&\sum_{i=1}^N\left|\lambda_i\left(-\Delta-\f{c_*}{|x|^2}-V_n\right)\right|^s	\\
=&\sum_{i=1}^M\left|\lambda_i\left(-\Delta-V_n\mathds{1}_{B_{R_n}}-\f{c_*}{|x|^2}\right)\right|^s+\sum_{i=1}^{N-M}\left|\lambda_i\left(-\Delta-V_n\mathds{1}_{\R^d\setminus B_{2R_n}}-\f{c_*}{|x|^2}\right)\right|^s+o(1)\\
=&\sum_{i=1}^M\left|\lambda_i\left(-\Delta-V_n\mathds{1}_{B_{R_n}}-\f{c_*}{|x|^2}\right)\right|^s+\sum_{i=1}^{N-M}\left|\lambda_i\left(-\Delta-V_n\mathds{1}_{\R^d\setminus B_{2R_n}}\right)\right|^s+o(1)\\
\leq&\sum_{i=1}^N\left|\lambda_i\left(-\Delta-V_n\mathds{1}_{B_{R_n}}-\f{c_*}{|x|^2}\right)\right|^s+\sum_{i=1}^{N}\left|\lambda_i\left(-\Delta-V_n\mathds{1}_{\R^d\setminus B_{2R_n}}\right)\right|^s+o(1)\\
\leq& C_{s,d}^{(N)}\int_{\R^d}\l(
V_n\mathds{1}_{B_{R_n}}\r)^{s+\f d2}dx+L_{s,d}^{(N)}\int_{\R^d}\l(
V_n\mathds{1}_{\R^d\setminus B_{2R_n}}\r)^{s+\f d2}dx+o(1)\\
=& C_{s,d}^{(N)}\int_{\R^d}\l(
V_n\mathds{1}_{B_{R_n}}\r)^{s+\f d2}dx+\big(C_{s,d}^{(N)}-\epsilon_0\big)\int_{\R^d}\l(
V_n\mathds{1}_{\R^d\setminus B_{2R_n}}\r)^{s+\f d2}dx+o(1)\\
=&C_{s,d}^{(N)}-\epsilon_0\int_{\R^d}\l(
V_n\mathds{1}_{\R^d\setminus B_{2R_n}}\r)^{s+\f d2}dx+o(1)\ \ \text{ as } \, n\to\infty,
\end{align*}
where the sequence $\{R_n\}$ is as in \eqref{a.4}.
This gives that
\begin{equation}\label{c.8}
C_{s,d}^{(N)}=\sum_{i=1}^N\left|\lambda_i\left(-\Delta-\f{c_*}{|x|^2}-V_n\right)\right|^s+\epsilon_0\int_{\R^d}\l(
V_n\mathds{1}_{\R^d\setminus B_{2R_n}}\r)^{s+\f d2}dx+o(1)\ \ \text{ as }\, n\to\infty.
\end{equation}
Note that $\{V_n\}$ is an optimizing sequence of the problem $C_{s,d}^{(N)}$ defined by \eqref{sup}. We thus deduce from \eqref{e0} and \eqref{c.8} that
	\begin{equation*}
		\lim_{n\to\infty}\int_{\R^d}\l(
		V_n\mathds{1}_{\R^d\setminus B_{2R_n}}\r)^{s+\f d2}dx=0.
	\end{equation*}
It further implies from \eqref{sup} and \eqref{a.4} that
\begin{equation}\label{2.47}
\al=\lim_{n\to\infty}\int_{\R^d}\l(
	V_n\mathds{1}_{B_{R_n}}\r)^{s+\f d2}dx=\lim_{n\to\infty}\int_{\R^d}
V_n^{{s+\f d2}}dx=1.
\end{equation}
We deduce from \eqref{2.43}, \eqref{a.4}  and \eqref{2.47} that
$$
\lim\limits_{R\to\infty}\liminf\limits_{n\to\infty}\int_{B_R}V_n^{s+\frac{d}{2}}dx=\lim\limits_{n\to\infty}\int_{B_{R_n}}V_n^{s+\frac{d}{2}}dx=\lim\limits_{n\to\infty}\int_{\R^d} V_n^{s+\frac{d}{2}}dx,
$$
which thus yields that \eqref{am} holds true. This completes  the proof of Proposition \ref{T2.1}.\qed
	
\section{Existence of Optimizers}\label{S3}
In this section, we shall complete the proof of Theorem \ref{exist}  on the existence of optimizers for the problem $C_{s,d}^{(N)}$ defined by \eqref{sup}.

To prove Theorem \ref{exist}, we first employ Proposition \ref{T2.1} to establish the following convergence of min-max levels.

\begin{Lemma}\label{eicon}
Let $\{V_n\geq0\}\subset L^{s+\f d2}(\R^d, \R)$ be a normalized optimizing sequence of the problem $C_{s,d}^{(N)}$ defined by \eqref{sup}, where $s>0$, $d\geq3$ and $N\in\N^+$. Then there exists a function $V\in L^{{s+\f d2}}(\R^d,\R)$ such that  up to a subsequence if necessary,
\begin{equation}\label{wc}
V_n\rightharpoonup V\text{ \ weakly  in \ }L^{{s+\f d2}}(\R^d,\R)\text{ \ as \ }n\to\infty,
\end{equation}
and
\begin{equation}\label{0.14}
\lim_{n\to\infty}\lambda_i\l(-\Delta-\f{c_*}{|x|^2}-V_n\r)= \lambda_i\l(-\Delta-\f{c_*}{|x|^2}-V\r),\ \, i=1,\cdots,N.
	\end{equation}
\end{Lemma}

\noindent{\bf Proof.}
Since the weak convergence \eqref{wc} follows immediately from the fact that \begin{equation}\label{v1}\|V_n\|_{s+\f d2}=1\text{ \ for all \ }n\in\N^+,
\end{equation}
the rest of the proof is to  prove \eqref{0.14}.

We first claim that 	
\begin{equation}\label{wsc}
	\limsup_{n\to\infty}\lambda_i\l(-\Delta-\f{c_*}{|x|^2}-V_n\r)\leq \lambda_i\l(-\Delta-\f{c_*}{|x|^2}-V\r),\ \, i=1,\cdots,N.
\end{equation}
For any fixed $i\in\{1,\cdots,N\}$, let  $u_1,\cdots,u_i\in   \mathcal Q$ form an  orthonormal family in $L^2(\R^d,\C)$. For any $n\in\N^+$, suppose $v_n:=\sum_{k=1}^ic_{kn}u_k$ is an optimizer of the following  finite-dimensional problem:
\begin{equation}\label{fp}
	\max_{\substack{u\in\operatorname{span}(u_1,\cdots,u_i)\\ \|u\|_2=1}}\int_{\R^d}\l(|\n u|^2-V_n|u|^2-\f{c_*}{|x|^2}|u|^2\r)dx,\,\ i\in\{1,\cdots,N\}.
\end{equation}
Since $\|v_n\|^2_2=\sum_{k=1}^i|c_{kn}|^2=1,$ there exists $c_k\in\C$ such that
\begin{equation*}
	c_{kn}\to c_k\ \text{ as \ }n\to\infty,\ \  k=1, \cdots, i.
\end{equation*}
We then deduce that
\begin{equation}\label{vsc} v_n\to v:=\sum_{k=1}^ic_{k}u_k \text{ \  strongly in\ }\mathcal Q \text{ \  as \ } n\to\infty.
\end{equation}
Applying H\"older's inequality and \eqref{gnh}, we thus deduce from \eqref{wc} and  \eqref{v1} that
\begin{align}\label{11}
%\begin{split}
&\Big|\int_{\R^d}V_n |v_n|^2dx-\int_{\R^d}V|v|^2dx\Big|\nonumber\\
\leq& \int_{\R^d}V_n\, \big||v_n|-|v|\big|\ \big(|v_n|+|v|\big)dx+\Big|\int_{\R^d}(V_n-V)|v|^2dx\Big|\\[1mm]
\leq&\|V_n\|_{s+\f d2}\cdot\big\||v_n|+|v|\big\|_{m}\cdot\big\||v_n|-|v|\big\|_{m}+o(1)\nonumber\\[1mm]
\leq &C\big\||v_n|+|v|\big\|_\mathcal Q^{\f{(m-2)d}{2m}}\cdot \big\||v_n|-|v|\big\|_\mathcal Q^{\f{(m-2)d}{2m}}+o(1)\nonumber\\[1mm]
\leq &C_1\big\|v_n-v\big\|_\mathcal Q^{\f{(m-2)d}{2m}}+o(1)\to0\ \ \text{ \ as \ }n\to\infty,\nonumber
%	\end{split}
\end{align}
where $m=(2d+4s)/(d-2+2s)$, and $C_1>0$ is independent of $n>0$. By the min-max principle,  we then calculate from \eqref{vsc} and \eqref{11} that
\begin{equation}\label{3.7A}
	\begin{split}
		&\limsup_{n\to\infty}\lambda_i\Big(-\Delta-V_n-\f{c_*}{|x|^2}\Big)\\
		\leq&\limsup_{n\to\infty}\max_{\substack{u\in\operatorname{span}(u_1,\cdots,u_i)\\ \|u\|_2=1}}
		\int_{\R^d}\Big(|\n u|^2-V_n|u|^2-\f{c_*}{|x|^2}|u|^2\Big)dx\\
		=&\limsup_{n\to\infty}
		\int_{\R^d}\Big(|\n v_n|^2-V_n|v_n|^2-\f{c_*}{|x|^2}|v_n|^2\Big)dx\\
		=&\int_{\R^d}\Big(|\n v|^2-V|v|^2-\f{c_*}{|x|^2}|v|^2\Big)dx\\
		\leq&\max_{\substack{u\in\operatorname{span}(u_1,\cdots,u_i)\\ \|u\|_2=1}}\int_{\R^d}\Big(|\n u|^2-V|u|^2-\f{c_*}{|x|^2}|u|^2\Big)dx,\ \ \ i\in\{1,\cdots,N\},
	\end{split}
\end{equation}
where the last inequality follows from the fact that $v=\sum_{k=1}^ic_{k}u_k\in\operatorname{span}(u_1,\cdots,u_i)$   satisfies $\|v\|_2=1$. Minimizing the right hand side of \eqref{3.7A} over $\{u_1,\cdots,u_i\}\subset   \mathcal Q$, we conclude that the claim \eqref{wsc} holds true.

Since \begin{equation*}\lim_{n\to\infty}{\sum_{i=1}^N\Big|\lambda_i\Big(-\Delta-\f{c_*}{|x|^2}-V_n\Big)\Big|^s}=C_{s,d}^{(N)}<\infty,
\end{equation*}
without
loss of generality, we may	suppose  that
\begin{equation}\label{a.5}
	\lim_{n\to\infty}\lambda_i\Big(-\Delta-\f{c_*}{|x|^2}-V_n\Big)=\Lambda_i\in(-\infty,0],\ \ i=1,\cdots,N.
\end{equation}
If $\Lambda_{i^0}=0$ holds for some $i^0\in\{1,2,\cdots,N\}$, then the weakly upper semi-continuity \eqref{wsc} immediately gives that \eqref{0.14} holds true for all $i=i^0,i^0+1,\cdots,N$.  Consequently, without loss of generality, we next suppose that $\Lambda_i<0$ holds for any $i=1,\cdots, N$.

 Choose $-E\in(\Lambda_N,\,0)$. It then follows from \eqref{a.5} that for sufficiently large $n>0$,
\begin{equation}\label{a.7} \lambda_{in}:=\lambda_i\l(-\Delta-\f{c_*}{|x|^2}-V_n\r)<-E<0,\ \ i=1,\cdots,N.
\end{equation}
Let $u_{in}$ be the corresponding $L^2$-normalized (real) eigenfunction of $ \lambda_{in}$, $i.e.$,
\begin{equation}\label{ee}
	\l(-\Delta-\f{c_*}{|x|^2}-V_n\r)u_{in}=\lambda_{in} u_{in}\text{ \ in \ }\R^d,\ \ i=1,\cdots,N.
\end{equation}
The same argument of \eqref{un} then yields that
\begin{equation}\label{3.12}
\text{the sequence }
\{u_{in}\}_n \text{\ is bounded uniformly in \ }  \mathcal Q,\ \ i=1,\cdots,N.
\end{equation}
Thus,  there exists $u_i\in\mathcal Q$ such that
\begin{equation}\label{wcq}
	u_{in}\rightharpoonup u_i \text{ \ weakly in $\mathcal Q\cap H^t(\R^d,\R)$ \, as\, $n\to\infty$},\ \ \forall\  t\in(0,1),\ \ i=1,\cdots ,N.
\end{equation}
Together with \eqref{ce},  we then conclude that
\begin{equation}\label{loccon}
	u_{in}\to u_{i}\text{\quad strongly in $L_{loc}^q(\R^d,\R)$ \,as\, $n\to\infty$},\ \ \forall \ q\in\Big[2,\f{2d}{d-2}\Big),\ \ i=1,\cdots ,N.
\end{equation}

Let $f_R(x)=f({x}/{R})$ for $R>0$, where $f(x)\in C_c^\infty (\R^d,[0,1])$ satisfies $f(x)=0$ for $|x|\leq1$ and $f(x)=1$ for $|x|\geq2$. Multiplying  \eqref{ee} by $f_R^2(x)u_{in}(x)$ and integrating over $\R^d$, we deduce that for $i=1, \cdots,N$,
\begin{equation}\label{a.8}
	\begin{split}
		\int_{\R^d}\Big[\n u_{in}\cdot\n(f_R^2{u}_{in})-\f{c_*}{|x|^2}|f_Ru_{in}|^2-\lambda_{in}|f_Ru_{in}|^2\Big]dx=\int_{\R^d}V_n|f_Ru_{in}|^2dx.
	\end{split}
\end{equation}
Note that
\begin{equation*}
	\begin{split}
		&\int_{\R^d}\Big[\n u_{in}\cdot\n(f_R^2{u}_{in})-\f{c_*}{|x|^2}|f_Ru_{in}|^2\Big]dx\\
		\geq&\int_{\R^d}f_R^2|\n{u}_{in}|^2dx-\f12 \int_{\R^d}|u_{in}|^2\Delta f_R^2dx-\f{C_1}{R^2}\\
		\geq&\int_{\R^d\setminus B_{2R}}|\n{u}_{in}|^2dx-\f{C_2}{R^2},\ \ i=1,\cdots,N,
	\end{split}
\end{equation*}
where $C_2>0$ is independent of $n>0$.	This then yields from \eqref{a.7} and \eqref{a.8} that for sufficiently large $n>0$,
\begin{equation}\label{M:a.9}
	\begin{split}
		\int_{\R^d\setminus B_{2R}}\big(|\n{u}_{in}|^2+E|u_{1,n}|^2\big)dx
		\leq C\|V_n\|_{L^{s+\f d2}(\R^d\setminus B_R)}+\f C{R^2},\ \ i=1,\cdots,N.
	\end{split}
\end{equation}
Note from Proposition \ref{T2.1} that the sequence $\{V_n\}$ is tight in the sense of \eqref{am}. It hence follows that
\begin{equation*}
    \lim_{R\to\infty}\lim_{n\to\infty}\|V_n\|_{L^{s+\f d2}(\R^d\setminus B_R)}=0.
\end{equation*}
We thus deduce from \eqref{M:a.9} that for all $i=1, \cdots, N$,
\begin{equation*}
\begin{split}
0\leq&\lim_{R\to\infty}\lim_{n\to\infty}\int_{\R^d\setminus B_{2R}}\big(|\n{u}_{in}|^2+|u_{1,n}|^2\big)dx\\
\leq&    \f{1}{\min\{1,E\}}\lim_{R\to\infty}\lim_{n\to\infty}\int_{\R^d\setminus B_{2R}}\big(|\n{u}_{in}|^2+E|u_{1,n}|^2\big)dx=0,
\end{split}
\end{equation*}
and hence for any fixed $q\in[2,\f{2d}{d-2})$, we deduce from \eqref{ce} and \eqref{3.12} that for $2\leq q<q'<\f{2d}{d-2}$,
\begin{equation*}
\begin{split}
0&\leq \lim_{R\to\infty}\lim_{n\to\infty}\|u_{in}\|_{L^q(B^c_{2R})}\\
&\leq \lim_{R\to\infty}\lim_{n\to\infty}\Big(\|u_{in}\|_{L^2(B^c_{2R})}^\theta\|u_{in}\|_{L^{q'}(B^c_{2R})}^{1-\theta}\Big)\\
&\leq \lim_{R\to\infty}\lim_{n\to\infty}\Big(\|u_{in}\|_{L^2(B^c_{2R})}^\theta\|u_{in}\|_{L^{q'}(\R^d)}^{1-\theta}\Big)\\
&\leq C \lim_{R\to\infty}\lim_{n\to\infty}\Big(\int_{\R^d\setminus B_{2R}}\big(|\n{u}_{in}|^2+|u_{1,n}|^2\big)dx\Big)^{\f{\theta}{2}}\\
&=0,
\end{split}
\end{equation*}
where we used interpolation inequality in the second inequality and $\theta$ is given by
$$\f{1}{q}=\f{\theta}{2}+\f{1-\theta}{q'}.$$
This further implies from \eqref{loccon} that
\begin{equation}\label{strcon}
	u_{in}\to u_i\ \text{ strongly in }L^q(\R^d,\R) \text{ as } n\to\infty,\ \  \forall \,2\leq q<\f{2d}{d-2}, \ \ i=1,\cdots,N.
\end{equation}
Consequently, we conclude that
\begin{equation}\label{ooo}
\langle  u_i,\,u_j\rangle_{2}=\lim_{n\to\infty}\langle u_{in},\,u_{jn}\rangle_{2}=\delta_{ij},\ \ i,\,j=1,\cdots,N.
\end{equation}

Applying H\"older's inequality, we obtain from \eqref{v1} and \eqref{strcon} that
\begin{equation}\label{3.16}
	\begin{split}
		&\Big|\int_{\R^d}V_n| u_{in}|^2dx-\int_{\R^d}V| u_{i}|^2dx\Big|\\ \leq&\|V_n\|_{s+\f d2}\cdot\big\||u_{in}|+|u_i|\big\|_{\f{2d+4s}{d-2+2s}}\cdot\big\||u_{in}|-|u_i|\big\|_{\f{2d+4s}{d-2+2s}}+\Big|\int_{\R^d}(V_n-V)|u_i|^2dx\Big|\\\to&0\text{ \ as \ }n\to\infty,\ i=1,\cdots,N,
	\end{split}
\end{equation}
where $s>0$, $d\geq3$ and $2<\f{2d+4s}{d-2+2s}<2^*=\f{2d}{d-2} $. By the min-max principle, we hence deduce from  \eqref{a.7}, \eqref{wcq} and \eqref{3.16}  that
\begin{equation*}
	\begin{split}
		\lim_{n\to\infty}\lambda_{1n}&=\lim_{n\to\infty}\int_{\R^d}\l(|\n u_{1n}|^2-V_n| u_{1n}|^2-\f{c_*}{|x|^2}| u_{1n}|^2\r)dx\\
		&\geq \int_{\R^d}\l(|\n u_1|^2-V| u_{1}|^2-\f{c_*}{|x|^2}| u_{1}|^2\r)dx\\
		&\geq \lambda_{1}\l(-\Delta-V-\f{c_*}{|x|^2}\r).
	\end{split}
\end{equation*}
Together with \eqref{wsc}, this yields from \eqref{a.7} that
\begin{equation}\label{3.18}
\lim_{n\to\infty}\lambda_{1n}=\lambda_1 \big(-\Delta-{c_*}|x|^{-2}-V\big).
\end{equation}
Furthermore, applying  the min-max principle again, one can verify  from \eqref{a.7}, \eqref{wcq},  \eqref{ooo} and \eqref{3.16} that
\begin{equation*}
	\begin{split}
		\lim_{n\to\infty}\sum_{i=1}^2\lambda_{in}&=\lim_{n\to\infty}\sum_{i=1}^2\int_{\R^d}\l(|\n u_{in}|^2-V_n| u_{in}|^2-\f{c_*}{|x|^2}| u_{2n}|^2\r)dx\\
		&\geq \sum_{i=1}^2\int_{\R^d}\l(|\n u_i|^2-V| u_{i}|^2-\f{c_*}{|x|^2}| u_{i}|^2\r)dx\\
		&\geq \sum_{i=1}^2\lambda_{i}\l(-\Delta-V-\f{c_*}{|x|^2}\r),
	\end{split}
\end{equation*}
which then implies from \eqref{3.18} that
\begin{equation*}
	\lim_{n\to\infty}\lambda_{2n}\geq\lambda_2\big(-\Delta-{c_*}|x|^{-2}-V\big).
\end{equation*}
Together with \eqref{wsc},  we thus conclude that
$$	\lim_{n\to\infty}\lambda_{2n}=\lambda_2\big(-\Delta-{c_*}|x|^{-2}-V\big).$$
By induction, we therefore  obtain that
\begin{equation*}
\lim_{n\to\infty}\lambda_{in}=\lambda_i\big(-\Delta-{c_*}|x|^{-2}-V\big),\  \ \ i=1,\cdots,N,
\end{equation*}
$i.e.$, \eqref{0.14} holds true. This completes the proof of Lemma \ref{eicon}.\qed

\begin{Remark}\label{R3.1}
In Lemma \ref{eicon},  we suppose  that $\{V_n\geq0\}$ is an optimizing sequence of \eqref{sup}. Actually, if  this assumption is replaced by
\begin{equation*}
\lim_{R\to\infty}\liminf_{n\to\infty}\int_{|x|\leq R}V_n^{s+\f d2}dx= \lim\limits_{n\to\infty}\int_{\R^d}V_n^{s+\f d2}dx,
	\end{equation*}
then Lemma \ref{eicon} still holds true.
\end{Remark}

Applying  Lemma \ref{eicon}, in the following we complete the proof of Theorem \ref{exist}.
\vspace{.10cm}

\noindent{\bf Proof of Theorem \ref{exist}.} Let the sequence $\{V_n\geq0\}$ be as in Lemma \ref{eicon}. By Lemma \ref{eicon}, it then follows from  \eqref{hlt} that for all $N\in\N^+$,
\begin{equation*}
	\begin{split}
		C_{s,d}^{(N)}&=\lim_{n\to\infty}\sum_{i=1}^N\l|\lambda_{i}\l(-\Delta-\f{c_*}{|x|^2}-V_n\r)\r|^s\\
		&=\sum_{i=1}^N\l|\lambda_{i}\l(-\Delta-V-\f{c_*}{|x|^2}\r)\r|^s\\
		&\leq C_{s,d}^{(N)}\int_{\R^d}V^{s+\f d2}dx\leq C_{s,d}^{(N)},
	\end{split}
\end{equation*}
due to the fact that $\|V\|_{s+\f d2}\leq\liminf_{n\to\infty}\|V_n\|_{s+\f d2}=1$, where the  function $V\in L^{{s+\f d2}}(\R^d,\R)$ is given by Lemma \ref{eicon}. This yields that $V$ is an optimizer of $C_{s,d}^{(N)}$, and the proof of  Theorem \ref{exist}  is therefore complete.\qed

\section{Analytical Properties of Optimizers }\label{S4}	
In this section, we shall prove Theorem \ref{analytic}  on the analytical properties of optimizers for the problem  \eqref{sup}.

As before, let $V_N$ be an optimizer of $C_{s,d}^{(N)}$ defined by \eqref{sup}, where $d\geq3$, $s>0$ and $N\in\N^+$. Suppose $M\in\N^+$ is the number of the negative eigenvalues (counted with multiplicity) for $-\Delta-\f{c_*}{|x|^2}-V_N$ in $\R^d$. If $M\geq N$, then it follows that $\lambda_N\leq\lambda_{M}<0$. If $M<N$, then $V_N$ is also an optimizer for $C_{s,d}^{(M)}$, and the gap $\lambda_M<\lambda_{M+1}=0$ holds true. Consequently, without loss of generality, in this section we always suppose that $\lambda_N<0$. We first prove the following lemma on the gap of the eigenvalues $\lambda_1,  \lambda_2, \cdots, \lambda_N$.

\begin{Lemma}\label{L3.3}
For any fixed $N\in\N^+$,	let $V_N$ be an $L^{s+\f d2}$-normalized optimizer of $C_{s,d}^{(N)}$ defined by \eqref{sup}, where $s>0$ and $d\geq3$. Suppose $\lambda_1\leq \lambda_2\leq\cdots\leq \lambda_N<0$ are the $N$-first negative eigenvalues of the operator $-\Delta-V_N-{c_*}{|x|^{-2}}$ in $\R^d$. Then we have	
\begin{equation}\label{gap}
\lambda_N<\lambda_{N+1}\leq0.
\end{equation}
\end{Lemma}

\noindent{\bf Proof.}
Similar to \cite[Theorem 11.8]{Lieb01},  one can verify  that the first eigenvalue $\lambda_1$ is simple. We thus deduce that if $N=1$, then $\lambda_{N}<\lambda_{N+1}$ holds true. Hence, we next consider the case $2\leq N\in \N^+$.

By contradiction, suppose $\lambda_{N}=\lambda_{N+1}<0$ holds for some $N\geq2$. Let $m\in \N^+$ be the multiplicity of $\lambda_{N}$,  and suppose $K\in\N^+$ is the smallest integer so that $\lambda_{N-K}<\lambda_N$. We then obtain that  $	 K<m$  and
\begin{equation}\label{3.22}
	\lambda_1<\lambda_2\leq\cdots\leq	\lambda_{N-K}<	 \lambda_{N-K+1}=\cdots=\lambda_N=\cdots=\lambda_{N-K+m}<\lambda_{N-K+m+1}.
\end{equation}
For convenience, we denote
\begin{equation}\label{HE}
	H:=-\Delta-\f{c_*}{|x|^2}-V_N\ \text{\ \ and\ }\ H_\epsilon:=-\Delta-\f{c_*}{|x|^2}-V_N-\epsilon h,
\end{equation}
where $\epsilon\in\R$ and $h\in L^{s+\f d2}(\R^d,\R)$ are arbitrary. It then follows from Remark \ref{R3.1} that \begin{equation}\label{ec0}
	\lim_{\epsilon\to0}\lambda_i(H_\epsilon)=\lambda_i<0,\ \  i=1,\cdots, N,
\end{equation}
where $\lambda_1<\lambda_2\leq\cdots\leq \lambda_N<0$ are the $N$-first negative eigenvalues of the operator $-\Delta-V_N-{c_*}{|x|^{-2}}$ in $\R^d$.

We now analyze the refined asymptotics of $\lambda_i(H_\epsilon)$ as $\epsilon\to0$, where $i=1,\,\cdots,\,N-K$. Choose negative numbers $a'<a<b<b'<0$ such that
\begin{equation*}
	a'<a<\lambda_1<\lambda_{N-K}<b<b'<\lambda_{N-K+1}<0.
\end{equation*}
We then obtain from \eqref{ec0} that  for sufficiently small $|\epsilon|>0$,
\begin{equation}\label{a.10}
\text{the operator}\ H_\epsilon\  \text{has exactly}\  N-K\  \text{negative eigenvalues in}\ [a,b].
\end{equation}
 Let $f\in C_c^{\infty}(\R^-,\R) $ be a  function supported in $[a',b']$ such that
\begin{equation}
\label{df}f(x)=(-x)^s,\ \ x\in[a,b],
\end{equation}
where $s>0$ is as \eqref{sup}. Applying the spectral mapping theorem (cf. \cite[Theorem \uppercase\expandafter{\romannumeral7}.1]{Sim1}),
we then deduce from \eqref{a.10} that
\begin{equation}\label{3.19}
	\operatorname{Tr}f(H_\epsilon )=\sum_{i=1}^{N-K}\big|\lambda_i(H_\epsilon)\big|^s \text{\  \ as \ }\epsilon\to0,
\end{equation}
where the operator $H_\epsilon$ is defined by \eqref{HE}. Moreover, we note from \eqref{HE} that
\begin{equation}\label{3.20}
\begin{split}
\operatorname{Tr}f(H_\epsilon)=\operatorname{Tr}f(H-\epsilon h)
		=\operatorname{Tr}f(H)-\epsilon \operatorname{Tr}\big( f'\left(H\right)h\big)+o(\epsilon)\text{ \ \ as \ }\epsilon\to0.
	\end{split}
\end{equation}
It then yields from \eqref{3.19} and \eqref{3.20} that
\begin{equation}\label{N-K}
	\begin{split}
		\sum_{i=1}^{N-K}\big|\lambda_i(H_\epsilon)\big|^s	&=\operatorname{Tr}f(H)-\epsilon\operatorname{Tr}\big( f'(H)h\big)+o(\epsilon)\\
		&=	\sum_{i=1}^{N-K}|\lambda_i|^s+\epsilon s\sum_{i=1}^{N-K}|\lambda_{i}|^{s-1}\int_{\R^d}|u_i|^2h\, dx+o(\epsilon)\text{ \ \ as \ }\epsilon\to0,
	\end{split}
\end{equation}
where $u_i$ is the corresponding eigenfunction of the eigenvalue $\lambda_i$ for $i=1,\cdots,N-K$.

%the eigenvalue $\lambda_i(H_\epsilon)$ has the asymptotic expansions

Let $\{u_{N-K+1},\cdots, u_{N-K+m}\}$ be an orthonormal basis of the corresponding  eigenspace for the eigenvalue $\lambda_N$, and suppose
\begin{equation}\label{3.23}	
	\mu_1^h\leq\cdots\leq\mu_m^h
\end{equation}
are the $m$ eigenvalues of the $m\times m$ Hermitian matrix
\begin{equation}\label{matrix1}
	\l(\int_{\R^d}\bar{u}_i u_jhdx\r)_{N-K+1\leq i,\,j\leq N-K+m}.
\end{equation}
For simplicity, we denote
\begin{equation}\label{3.32}
\hat{\mu}_1^h:=\mu_m^h,\ \ \hat{\mu}_2^h:=\mu_{m-1}^h,\ \ \cdots,\ \ \hat{\mu}_m^h:=\mu_1^h.
\end{equation}
By  \cite[Chapter 8, Theorem 2.6]{Kat95},  we then get from \eqref{3.22} and \eqref{ec0} that for $i=N-K+1,\cdots,N-K+m$,
\begin{equation*}\label{ae1}
	\lambda_i(H_\epsilon)=\lambda_N-\epsilon \mu_{i-N+K}^h+o(\epsilon) \ \ \text{ as }\ \epsilon\to0^-,
\end{equation*}
and
\begin{equation*}\label{ae2}
	\lambda_i(H_\epsilon)=\lambda_N-\epsilon \hat{\mu}_{i-N+K}^h+o(\epsilon) \ \ \text{ as }\ \epsilon\to0^+,
\end{equation*}
which further imply that
\begin{equation}\label{N-K+11}
\begin{split}
\sum_{i=N-K+1}^{N}\big|\lambda_i(H_\epsilon)\big|^s&=\sum_{i=N-K+1}^{N}\big|\lambda_N-\epsilon{\mu}_{i-N+K}^h+o(\epsilon)\big|^s\\
&=\sum_{i=N-K+1}^{N}\Big[\big|\lambda_N\big|^s+\epsilon s|\lambda_{N}|^{s-1}{\mu}_{i-N+K}^h+o(\epsilon)\Big]\\
&=\sum_{i=N-K+1}^{N}\big|\lambda_N\big|^s+\epsilon s|\lambda_{N}|^{s-1}\sum_{k=1}^{K} \mu_{k}^h+o(\epsilon)\ \ \text{ as }\ \epsilon\to0^-,
\end{split}
\end{equation}
and similarly,
\begin{equation}\label{N-K+12}
	\begin{split}
		\sum_{i=N-K+1}^{N}\big|\lambda_i(H_\epsilon)\big|^s&=\sum_{i=N-K+1}^{N}\big|\lambda_N\big|^s+\epsilon s|\lambda_{N}|^{s-1}\sum_{k=1}^{K} \hat{\mu}_{k}^h+o(\epsilon)\ \ \text{ as }\ \epsilon\to0^+.
	\end{split}
\end{equation}
We then derive from \eqref{sup} and \eqref{N-K} that
\begin{equation}\label{3.131}
	\begin{split}
	\sum_{i=1}^{N}\big|\lambda_i(H_\epsilon)\big|^s
	=&C_{s,d}^{(N)}+\epsilon s\sum_{i=1}^{N-K}|\lambda_{i}|^{s-1}\int_{\R^d}|u_i|^2hdx\\
	&+\epsilon s|\lambda_{N}|^{s-1}\sum_{k=1}^{K} \mu_{k}^h+o(\epsilon)\ \text{ \ as \ }\epsilon\to0^-,
	\end{split}\end{equation}
and
\begin{equation}\label{3.132}
	\begin{split}
	\sum_{i=1}^{N}\big|\lambda_i(H_\epsilon)\big|^s
	=&C_{s,d}^{(N)}+\epsilon s\sum_{i=1}^{N-K}|\lambda_{i}|^{s-1}\int_{\R^d}|u_i|^2hdx\\
	&+\epsilon s|\lambda_{N}|^{s-1}\sum_{k=1}^{K} \hat{\mu}_{k}^h+o(\epsilon)\ \text{ \ as \ }\epsilon\to0^+.
	\end{split}\end{equation}

On the other hand, we deduce from  \eqref{hlt} that
\begin{equation}\label{optimal}
\begin{split}
\sum_{i=1}^N\big|\lambda_i(H_\epsilon)\big|^s
\leq &C_{s,d}^{(N)}{\int_{\R^d}\left(V_N+\epsilon h\right)_+^{s+\f d2}dx}\\
=&C_{s,d}^{(N)}+\epsilon C_{s,d}^{(N)}\f{d+2s}{2}\int_{\R^d}V_N^{s+\f d2-1}hdx+o(\epsilon)\ \text{ \ as \ }\epsilon\to0.
\end{split}\end{equation}
Applying \eqref{3.131} and \eqref{3.132}, it yields from \eqref{optimal} that
\begin{equation}\label{3.381}
\begin{split}
&s|\lambda_{N}|^{s-1}\sum_{k=1}^{K} \mu_{k}^h\\
\geq&C_{s,d}^{(N)}\f{d+2s}{2}\int_{\R^d}V_N^{s+\f d2-1}hdx-s\sum_{i=1}^{N-K}|\lambda_{i}|^{s-1}\int_{\R^d}|u_i|^2hdx+o(1)\ \text{ \ as \ }\epsilon\to0^-,
\end{split}
\end{equation}
and
\begin{equation}\label{3.39}
\begin{split}
&s|\lambda_{N}|^{s-1}\sum_{k=1}^{K} \hat\mu_{k}^h\\
\leq&C_{s,d}^{(N)}\f{d+2s}{2}\int_{\R^d}V_N^{s+\f d2-1}hdx-s\sum_{i=1}^{N-K}|\lambda_{i}|^{s-1}\int_{\R^d}|u_i|^2hdx+o(1)\ \text{ \ as \ }\epsilon\to0^+.
\end{split}
\end{equation}
It then follows from \eqref{3.381} and \eqref{3.39} that
\begin{equation}\label{3.38}
\sum_{k=1}^{K} \hat{\mu}_{k}^h\leq\sum_{k=1}^{K} {\mu}_{k}^h.
\end{equation}
Since $K<m$, we deduce from \eqref{3.23}, \eqref{3.32} and  \eqref{3.38} that
$$\mu_1^h=\cdots=\mu_m^h.$$
This shows  that the Hermitian matrix matrix \eqref{matrix1} is  a multiplicity of the identity matrix, and hence
\begin{equation*}
	\int_{\R^d}\bar{u}_i u_jhdx=0,\ \ \int_{\R^d}\big(|u_i|^2-|u_j|^2\big)hdx=0,\ \  \forall \ 1\leq i\neq j\leq m, \ h\in L^{s+\f d2}(\R^d,\R).
\end{equation*}
Since $ h\in L^{s+\f d2}(\R^d,\R)$ is arbitrary, we obtain  that
\begin{equation*}
\bar{u}_i(x) u_j(x)\equiv0\ \text{ and }\  |u_i(x)|^2\equiv|u_j(x)|^2\ \ \mbox{in}\ \, \R^d, \ \  \forall\  1\leq i\neq j\leq m,
\end{equation*}
which further yields that
$$u_i(x)\equiv0\ \  \,\text{in}\ \ \R^d, \ \  i=1,\cdots,m.$$
This is however impossible, because  $u_i$ is an eigenfunction corresponding to the eigenvalue $\lambda_i$, where $i=1,\cdots,m$. Thus, the assumption $\lambda_{N}=\lambda_{N+1}$ cannot occur for any $N\geq2$, and the proof of Lemma \ref{L3.3} is therefore complete.	\qed

\vspace{.15cm}
Following Lemma \ref{L3.3}, we next finish the proof of Theorem \ref{analytic}.
\vspace{.15cm}

\noindent{\bf Proof of Theorem \ref{analytic}.} We first note that the estimate \eqref{Gap} follows directly from  Lemma \ref{L3.3}.

Under the assumptions of
Theorem \ref{analytic}, let
\begin{equation}\label{3.42}
(\lambda_i, u_i) \text{ \ be the $i$-th eigenpair of  the operator $-\Delta-V_N-{c_*}{|x|^{-2}}$ in $\R^d$},
\end{equation}
where $i=1, \cdots,N$.
Since $\lambda_N<\lambda_{N+1}(-\Delta-V_N-{c_*}{|x|^{-2}})\leq 0$, we can choose negative numbers $a'<a<b<b'$, so that
\begin{equation*}
	a'<a<\lambda_1<\lambda_N<b<b'<\lambda_{N+1}\leq0.
\end{equation*}
The similar argument of \eqref{N-K} then gives that% for  sufficiently small $|\epsilon|>0$,
\begin{equation}\label{330}
	\begin{split}
		\sum_{i=1}^{N}\big|\lambda_i(H_\epsilon)\big|^s = \sum_{i=1}^{N}|\lambda_i|^s+\epsilon s\sum_{i=1}^{N}|\lambda_{i}|^{s-1}\int_{\R^d}|u_i|^2h\, dx+o(\epsilon)\text{ \ \ as \ }\epsilon\to0,
	\end{split}
\end{equation}
where the operator $H_\epsilon$ is as in \eqref{HE}, and the function $h(x)\in L^{s+\f d2}(\R^d,\R)$ is arbitrary.
Substituting \eqref{330} into \eqref{optimal},  one gets that for  sufficiently small $|\epsilon|>0$,
\begin{equation*}
	\epsilon\l(s\sum_{i=1}^N|\lambda_{i}|^{s-1}\int_{\R^d}|u_i|^2hdx-C_{s,d}^{(N)}\f{d+2s}{2}\int_{\R^d}V_N^{s+\f d2-1}hdx\r)+o(\epsilon)\leq 0.
\end{equation*}
Since the above inequality holds for both positive and negative  $\epsilon$, we deduce that
\begin{equation*}
	s\sum_{i=1}^N|\lambda_{i}|^{s-1}\int_{\R^d}|u_i|^2hdx=C_{s,d}^{(N)}\f{d+2s}{2}\int_{\R^d}V_N^{s+\f d2-1}hdx.
\end{equation*}
Because $h\in L^{s+\f d2}(\R^d,\R)$ is arbitrary, this yields that
\begin{equation}\label{eV}
	V_N=\Big(\f{2s}{(d+2s)C_{s,d}^{(N)}}\sum_{i=1}^N|\lambda_{i}|^{s-1}|u_i|^2\Big)^{\f{2}{2s+d-2}},
\end{equation}
and hence both \eqref{ev} and \eqref{A:Gap} hold true in view of \eqref{3.42}.

We next prove that
$V_N\in C^2(\R^d\setminus\{0\}),$
for which it suffices to show that
$$V\in C^2 \big(B_a(x_0)\big)\ \text{ for any fixed}\ x_0\in\R^d\setminus\{0\}, \  \  \text{where} \ a:=\frac{|x_0|}{2}>0.
$$
For convenience, the dual exponent of $s+\f d2$ is denoted by
\begin{equation}
	p:=\l(s+\f d2\r)'=\f{2s+d}{2s+d-2}.
\end{equation}
Since $V_N\in L^{s+\f d2}(\R^d)$, we deduce from \eqref{eV} that
$u_i\in L^{(d+2s)(p-1)}(B_a(x_0))$.
Suppose that $u_i\in L^r(B_a(x_0))$ holds for some $r\geq (d+2s)(p-1)$, where  $i=1, \cdots N$. It then follows that $V\in L^{\f{r}{2(p-1)}}(B_a(x_0))$, and thus $V_Nu_i\in L^{\f{r}{2p-1}}(B_a(x_0))$, where $i=1,\cdots,N$. Recall that  $(\lambda_i, u_i)$ is an eigenpair of  the operator $-\Delta-V_N-{c_*}{|x|^{-2}}$ in $\R^d$, namely,
\begin{equation}\label{eigene}
	-\Delta u_i-\l(\f{c_*}{|x|^2}+\lambda_i\r)u_i=V_Nu_i\ \ \ \text{in}\ \, \R^d,\ \   i=1,\cdots,N,
\end{equation}
where  $\lambda_1<\lambda_2\leq\cdots\leq\lambda_N<0$. By the $L^p$ theory (cf.\,\cite{GT}), we then deduce from \eqref{eigene} that $u_i\in W^{2,\f{r}{2p-1}}_{loc}\big(B_a(x_0)\big)$ holds for $i=1,\cdots,N$. By Sobolev's embedding  theorem (cf.\,\cite[Theorem 4.12]{Adams}), this yields that
$u_i\in L^q_{loc}\big(B_a(x_0)\big)$ holds for $i=1,\cdots,N$, where $q\in[r,\infty)$ holds for $r\geq\f d2(2p-1)$, and
$q\in[r,r^*]$ satisfying $r^*:=r\big(\f{d}{d(2p-1)-2r}\big)$ holds for $r<\f d2(2p-1)$. Since $t_0=d(p-1)$ is an unstable fixed point of the function
$$f(t)=\f{dt}{d(2p-1)-2t},$$
we deduce that
$$r^*>r, \ \ r\geq (d+2s)(p-1)>d(p-1)=t_0.$$
Repeating the above argument at  finite  times, it then follows that $u_i\in L^q_{loc}(B_a(x_0))$ holds for all $q<\infty$, and hence $V_Nu_i\in L^q_{loc}(B_a(x_0))$ holds for any $q<\infty$, where $i=1,\cdots,N$. Using the $L^p$ theory again, we obtain from \eqref{eigene} that $u_i\in W^{2,q}_{loc}(B_a(x_0))$ holds for all $q< \infty$. By Sobolev's embedding theorem,  this implies that $u_i\in C_{loc}^{1,\al}(B_a(x_0))$ holds for some $0<\al<1$. In particular, we have $V_N\,u_i\in C^\al_{loc}(B_a(x_0))$. By the Schauder estimate (cf. \cite[Theorem 6.2]{GT}), we thus conclude from \eqref{eigene} that $u_i\in C^{2,\al}_{loc}(B_a(x_0))$ holds for $i=1,\cdots,N$. Since $|u_1|^2>0$, we also deduce from \eqref{eigene} with $i=1$ that $V_N\in C^2(B_a(x_0))$. This further proves $V_N\in C^2(\R^d\setminus\{0\})$.

We finally establish the  exponential decay \eqref{1.20}. Applying Kato's inequality (cf. \cite[Theorem X.27]{Sim2}), we deduce from \eqref{eigene} that
\begin{equation}\label{U-i}
	\left(-\Delta-V_N\right)|u_i|\leq\Big(\lambda_i+\f{c_*}{|x|^2}\Big) |u_i| \ \  \text{ in}\  \, \R^d,\ \  i=1, \cdots, N.
\end{equation}
Since $\lambda_i<0$, this implies that  there exists a sufficiently large constant $R_1>0$ such that
\begin{equation}\label{3.46}
\left(-\Delta-V_N\right)|u_i|\leq0\ \  \text{ in}\  \, B_{R_1}^c,\ \  i=1,\cdots,N.
\end{equation}
Utilizing De Giorgi-Nash-Moser theory  (cf. \cite[Theorem 4.1]{HL}), since $V_N\in L^{s+\f d2}(\R^d,\R)$, we obtain from \eqref{3.46}  that
\begin{equation*}
	\|u_i\|_{L^\infty(B_1(y))}\leq C\|u_i\|_{L^2(B_2(y))}\ \  \ \text{for any}\ B_2(y)\subset B_{R_1}^c,\ \   i=1,\cdots,N,
\end{equation*}
where $C>0$ depends only on $d\geq3$ and $\|V_N\|_{s+\frac{d}{2}}$.
Together with $u_i\in L^2(\R^d,\C)$, we deduce that
\begin{equation*}
\lim\limits_{|x|\to\infty}u_i(x)=0,\ \  i=1, \cdots, N.
\end{equation*}
This further yields from \eqref{eV} and \eqref{U-i} that there exists a large constant  $R_2$ satisfying $R_2>R_1>0$ such that
\begin{equation}\label{3.48}
\Big(-\Delta-\frac{3}{4}|\lambda_i|\Big)|u_i|\leq0\ \ \, \text{in}\ \, B_{R_2}^c,\ \  i=1,\cdots,N.
\end{equation}
Using the standard comparison  principle, we hence obtain from \eqref{3.48} that there exists a sufficiently large constant $R>0$ such that
\begin{equation*}\label{uved}
	|u_i(x)|\leq C\e ^{\f{-\sqrt{|\lambda_i|}}{2}|x|}\ \,\ \text{in}\  \, B_R^c,\  \ i=1, \cdots, N.
\end{equation*}
We thus conclude \eqref{1.20} in view of \eqref{eV}. This completes the proof of Theorem \ref{analytic}.\qed

\vspace{7pt}

\noindent{\bf Acknowledgements:} The authors thank Prof. Phan Th\'anh Nam and Dr. Long
 Meng very much for their fruitful discussions on the present paper.


\vspace{.37cm}



\begin{thebibliography}{45}
		
\bibitem{Adams} R.~A. Adams and J.~J.~F. Fournier,  Sobolev Spaces,  Pure and Applied Mathematics  Vol.~140, Elsevier/Academic Press, Amsterdam, 2003.
		
\bibitem{Lieb94} V.~Bach, E.~H. Lieb, M.~Loss and J.~P. Solovej, {\it There are no unfilled shells in unrestricted Hartree-Fock theory}, Phys. Rev. Lett. {\bf 72} (1994),  2981--2983.
		
\bibitem{BRS} G. M. Bisci, V.~D. Radulescu and R. Servadei,  Variational Methods for Nonlocal Fractional Problems, Encyclopedia of Mathematics and its Applications 162, Cambridge Univ. Press, Cambridge, 2016.
		
		
%\bibitem{Bre11} H.~Brezis,  Functional Analysis, Sobolev Spaces and Partial	Differential Equations, Universitext, Springer, New York, 2011.
		
		
\bibitem{CFKS87} H.~L. Cycon, R.~G. Froese, W.~Kirsch and B.~Simon,  Schr\"odinger Operators with Application to Quantum Mechanics and Global Geometry, Texts and Monographs in Physics, Springer-Verlag, Berlin, 1987.
		
		
\bibitem{Frank24} G.~K. Duong, R.~L. Frank, T.~M.~T. Le, P.~T. Nam and P. T. Nguyen, {\it Cwikel-Lieb-Rozenblum type inequalities for Hardy-Schr\"odinger operator}, J. Math. Pures Appl. (9) {\bf 190} (2024),  103598.
		
		
\bibitem{EF06} T.~Ekholm and R.~L. Frank, {\it On {L}ieb-Thirring inequalities for Schr\"odinger operators with virtual level}, Comm. Math. Phys. {\bf 264} (2006), 725--740.
		
\bibitem{Evans} L. C. Evans, Partial Differential Equations, second edition,  Graduate Studies in Mathematics, 19, Amer. Math. Soc. Providence, RI, 2010.		
		
\bibitem{Frank09} R.~L. Frank, {\it A simple proof of {H}ardy-{L}ieb-{T}hirring 	inequalities}, Comm. Math. Phys. {\bf 290} (2009),   789--800.
		
\bibitem{Frank} R.~L. Frank, {\it The Lieb-Thirring inequalities: recent results and open problems, Nine mathematical challenges---an elucidation},  Proc. Sympos. Pure Math. 104, AMS, Providence, RI, (2021), 45--86.
		
\bibitem{FGL21} R.~L. Frank, D.~Gontier and M.~Lewin, {\it The nonlinear {S}chr\"odinger equation for orthonormal functions {II}: {A}pplication to {L}ieb-{T}hirring inequalities}, Comm. Math. Phys. {\bf384} (2021),   1783--1828.
		
		
\bibitem{FGL24} R. L. {Frank}, D. {Gontier} and M.~{Lewin}, {\it Optimizers for the finite-rank Lieb-Thirring inequality}, Amer. J. Math.   (2024), to appear.
		
		
\bibitem{FLS08} R. L. Frank, E. H. Lieb and R. Seiringer, {\it Hardy-{L}ieb-{T}hirring inequalities for fractional {S}chr\"odinger operators}, J. Amer. Math. Soc. {\bf 21} (2008),  ~925--950.
		

\bibitem{GT} D.~Gilbarg and N.~S. Trudinger, Elliptic Partial Differential Equations of Second Order, Classics in Math., Springer-Verlag, Berlin, 2001. 		
		
		
\bibitem{GLN21} D.~Gontier, M.~Lewin and F.~Q. Nazar, {\it The nonlinear {S}chr\"odinger 	equation for orthonormal functions: existence of ground states}, Arch. Ration. Mech. Anal. {\bf240} (2021),  1203--1254.
		
		
\bibitem{HL} Q.~Han and F.~Lin,  Elliptic Partial Differential Equations, Courant Lecture Notes in Math. Vol.~1, CIMS, New York; AMS, Providence, RI, 2011.
		
		
%\textcolor{red}{\bibitem{HI77} I.~W. Herbst, {\it Spectral theory of the operator {$(p^2+m^2)^{1/2}-Ze^2/r$}}, Comm. Math. Phys. {\bf53} (1977),  ~285--294.}
				
\bibitem{Kat95} T.~Kato,   Perturbation Theory for Linear Operators, Classics in	Math., Springer-Verlag, Berlin, 1995.
		
		
\bibitem{KMK} M.~K. Kwong, {\it Uniqueness of positive solutions of $\Delta u-u+u^p=0$ in ${\R}^n$}, Arch. Ration. Mech. Anal. {\bf 105} (1989), no.~3, 243--266.
		
		
%\bibitem{Lewin11} M.~Lewin, {\it Geometric methods for nonlinear many-body quantum systems}, J. Funct. Anal. {\bf 260} (2011),  3535--3595.
		
		%\bibitem{Lieb81}{ E.~H. Lieb}, {\it Variational principle for many-fermion systems}, Phys. Rev. Lett. {\bf46} (1981),  ~457--459.
		
		
\bibitem{Lieb01} E.~H. Lieb and M.~Loss,   Analysis,  Graduate Studies in	Math. Vol. 14, AMS, Providence, RI, 2001.
		
		
\bibitem{LS10} E. H. Lieb and R. Seiringer,  The Stability of Matter in Qantum Mechanics, Cambridge University Press, Cambridge, 2010.
		
		
\bibitem{LT75} E.~H. Lieb and W.~E. Thirring, {\it Bound for the kinetic energy of fermions which proves the stability of matter}, Phys. Rev. Lett. {\bf 35} (1975), 687--689.
		
		
\bibitem{LT76}  E.~H. Lieb and W.~E. Thirring, {\it Inequalities for the moments of the eigenvalues of the schrodinger hamiltonian and their relation to sobolev inequalities}, Studies in Mathematical Physics,  Princeton University Press, Princeton (1976),  269--304.
		
		
\bibitem{Lions841} P. L. Lions, {\it The concentration-compactness principle in the calculus 	of variations. The locally compact case. I}, Ann. Inst. H. Poincar\'e Anal. Non Lin\'eaire, {\bf 1} (1984),  109--145.
		
%\textcolor{red}{\bibitem{Nam16} D.~Lundholm, P.~T. Nam and F.~Portmann, {\it Fractional {H}ardy-{L}ieb-{T}hirring and related inequalities for interacting systems}, Arch. Ration. Mech. Anal. {\bf 219} (2016),  1343--1382.}
		
		%\bibitem{Mor79}{ J.~D. III. Morgan}, {\it Schr\"odinger operators whose potentials have separated singularities}, J. Operator Theory, {\bf1} (1979),  ~109--115.
		
		%\bibitem{MS80}{ J.~D. III. Morgan and B.~Simon}, {\it Behavior of molecular potential energy curves for large nuclear separations}, International Journal of Quantum Chemistry, {\bf17} (1980),  ~1143--1166.
		
\bibitem{Nam21} D. Mukherjee, P. T. Nam and P. T. Nguyen, {\it Uniqueness of ground state and minimal-mass blow-up solutions for focusing NLS with Hardy potential}, J. Funct. Anal. {\bf 281} (2021),    109092.
		
\bibitem{Sim1} M. C. Reed and B. Simon, Methods of Modern Mathematical Physics. I. Functional Analysis, Academic Press, New York-London, 1972.
		
\bibitem{Sim2} M. C. Reed and B. Simon,   { Methods of Modern Mathematical Physics. II. Fourier Analysis, Self-Adjointness}, Academic Press, New York-London, 1975.
		
		
\bibitem{Sim05} B.~Simon,  Trace Ideals and Their Applications, Math. Surveys and Monographs Vol.~120, AMS, 	Providence, RI,   2005.
		
		%\bibitem{SJP91} { J.~P. Solovej}, {\it Proof of the ionization conjecture in a reduced {H}artree-{F}ock model}, Invent. Math. {\bf104} (1991),  ~291--311.

\bibitem{TGG}  G.~G. Talenti, {\it Best constant in Sobolev inequality}, Ann. Mat. Pura Appl. (4) {\bf 110} (1976), 353--372.

		
\bibitem{Ti99} T.~Weidl, {\it Remarks on virtual bound states for semi-bounded operators}, Comm. Partial Differential Equations {\bf 24} (1999),  25--60.
		
\end{thebibliography}
\end{document}